# MEAN-FIELD BACKWARD STOCHASTIC DIFFERENTIAL EQUATIONS: A LIMIT APPROACH

By Rainer Buckdahn, Boualem Djehiche, Juan Li[1,2]
and Shige Peng[2]

*Université de Bretagne Occidentale, Royal Institute of Technology,
Shandong University at Weihai and Shandong University*

Mathematical mean-field approaches play an important role in different fields of Physics and Chemistry, but have found in recent works also their application in Economics, Finance and Game Theory. The objective of our paper is to investigate a special mean-field problem in a purely stochastic approach: for the solution $(Y, Z)$ of a mean-field backward stochastic differential equation driven by a forward stochastic differential of McKean–Vlasov type with solution $X$ we study a special approximation by the solution $(X^N, Y^N, Z^N)$ of some decoupled forward–backward equation which coefficients are governed by $N$ independent copies of $(X^N, Y^N, Z^N)$. We show that the convergence speed of this approximation is of order $1/\sqrt{N}$. Moreover, our special choice of the approximation allows to characterize the limit behavior of $\sqrt{N}(X^N - X, Y^N - Y, Z^N - Z)$. We prove that this triplet converges in law to the solution of some forward–backward stochastic differential equation of mean-field type, which is not only governed by a Brownian motion but also by an independent Gaussian field.

**1. Introduction.** Our present work on a stochastic limit approach to a mean-field problem is motivated on the one hand by classical mean-field approaches in Statistical Mechanics and Physics (e.g., the derivation of Boltzmann or Vlasov equations in the kinetic gas theory), by similar approaches in Quantum Mechanics and Quantum Chemistry (in particular the density

Received February 2008; revised October 2008.

[1]Supported in part by the NSF of P. R. China (No. 10701050; 10671112), Shandong Province (No. Q2007A04).

[2]National Basic Research Program of China (973 Program) (No. 2007CB814904; 2007CB814906).

functional models or Hartree and Hartree–Fock type models) but also by a recent series of papers by Lasry and Lions (see [7] and the references therein). On the other hand it stems its motivation from partial differential equations of McKean–Vlasov type, which have found a great interest in recent years and have been studied with the help of stochastic methods by several authors.

Lasry and Lions introduced in recent papers a general mathematical modeling approach for high-dimensional systems of evolution equations corresponding to a large number of "agents" (or "particles"). They extended the field of such mean-field approaches also to problems in Economics and Finance as well as to the game theory: here they studied $N$-players stochastic differential games [7] and the related problem of the existence of Nash equilibrium points, and by letting $N$ tend to infinity they derived in a periodic setting the mean-field limit equation.

On the other hand, in the last years models of large stochastic particle systems with mean-field interaction have been studied by many authors; they have described them by characterizing their asymptotic behavior when the size of the system becomes very large. The reader is referred, for example, to the works by Bossy [1], Bossy and Talay [2], Chan [4], Kotelenez [6], Méléard [8], Overbeck [10], Pra and Hollander [14], Sznitman [15, 16], Talay and Vaillant [17] and all the references therein. They have shown that probabilistic methods allow to study the solution of the linear McKean–Vlasov PDE (see, e.g., Méléard [8]). With the objective to give a stochastic interpretation to such semilinear PDEs we introduced in [3] the notion of a mean-field backward stochastic differential equation: backward stochastic differential equations (BSDEs) have been introduced by Pardoux and Peng in their pioneering papers [11] and [12], in which they proved in particular that, in a Markovian framework, the solution of a BSDE describes the viscosity solution of the associated semilinear PDE. The reader interested in more details is also referred to the paper by El Karoui, Peng and Quenez [5] and the references therein. However, in order to generalize this stochastic interpretation to semilinear McKean–Vlasov PDEs we have had to study a new type of BSDE, which takes into account the specific, nonlocal structure of these PDEs; we have called this new type of backward equations mean-field BSDE (MFBSDE).

The objective of the present paper is to investigate a special mean-field problem in a purely stochastic approach: for the solution $(Y, Z)$ of a MF-BSDE driven by a forward SDE of McKean–Vlasov type with solution $X$ we study a special approximation by the solution $(X^N, Y^N, Z^N)$ of some decoupled forward–backward equation which coefficients are governed by $N$ independent copies of $(X^N, Y^N, Z^N)$. We show that the convergence speed of this approximation is of order $1/\sqrt{N}$. Moreover, our special choice of the approximation allows to characterize the limit behavior of $\sqrt{N}(X^N -$ $X, Y^N - Y, Z^N - Z$). We prove that this triplet converges in law to the solution of some forward–backward SDE of mean-field type, which is not only governed by a Brownian motion but also by an independent Gaussian field.

To be more precise, for a given finite time horizon $T > 0$, a $d$-dimensional Brownian motion $W = (W_t)_{t \in [0,T]}$ and a driving $d$-dimensional adapted stochastic process $X = (X_t)_{t \in [0,T]}$ we consider the BSDE of mean-field type

$$
\begin{aligned}
(1.1) \qquad dY_t &= -E[f(t, \lambda, \Lambda_t)]_{|\lambda = \Lambda_t} \, dt + Z_t \, dW_t, \qquad t \in [0, T], \\
Y_T &= E[\Phi(x, X_T)]_{x = X_T}
\end{aligned}
$$

with $\Lambda_t = (X_t, Y_t, Z_t)$, where $(Y, Z) = (Y_t, Z_t)_{t \in [0,T]}$ denotes the solution of the above equation. This equation also can be written in the form

$$
\begin{aligned}
dY_t &= -\int_{\mathbb{R}^d \times \mathbb{R} \times \mathbb{R}^d} f(t, \Lambda_t, \lambda) P_{\Lambda_t}(d\lambda) \, dt + Z_t \, dW_t, \qquad t \in [0, T], \\
Y_T &= \int_{\mathbb{R}^d} \Phi(X_T, x) P_{X_T}(dx).
\end{aligned}
$$

Such type of BSDEs has been studied recently by Buckdahn, Li and Peng [3] in a more general framework. In [3] for the "Markovian-like" case in which the process $X$ is the solution of a forward equation of McKean–Vlasov type (of course, $X$ is not a Markov process) it has been shown that such a BSDE gives a stochastic interpretation to the associated nonlocal partial differential equations. In the present work, given an arbitrary sequence of adapted processes $X^N = (X_t^N)_{t \in [0,T]}$ such that $E[\sup_{t \in [0,T]} |X_t^N - X_t|^2] \to 0$, we study the approximation of the above MFBSDE by a backward equation of the form

$$
\begin{aligned}
(1.2) \qquad dY_t^N &= -\frac{1}{N} \sum_{j=1}^{N} f(t, \Lambda_t^N, \Lambda_t^{N,j}) \, dt + Z_t^N \, dW_t, \qquad t \in [0, T], \\
Y_T^N &= \frac{1}{N} \sum_{j=1}^{N} \Phi(X_T^N, X_T^{N,j}).
\end{aligned}
$$

This approximating BSDE is driven by the i.i.d. sequence $\Lambda^{N,j} = (X^{N,j}, Y^{N,j}, Z^{N,j})$, $1 \leq j \leq N$, of triplets of stochastic processes, which are supposed to obey the same law as $\Lambda^N = (X^N, Y^N, Z^N)$ and to be independent of the Brownian motion $W$. In Section 3 we will introduce a natural framework in which such BSDEs can be solved easily, and we will show that the solution $(Y^N, Z^N)$ converges to that of (1.1) (Theorem 3.1). For the study of the speed of the convergence we assume that the process $X$ is the solution of the stochastic differential equation (SDE) of McKean–Vlasov type

$$
\begin{aligned}
(1.3) \qquad dX_t &= E[\sigma(t, x, X_t)]_{|x = X_t} \, dW_t + E[b(t, x, X_t)]_{|x = X_t} \, dt, \\
X_0 &= x_0 \in \mathbb{R}^m, \qquad t \in [0, T],
\end{aligned}
$$

and we approximate $X$ in the same spirit as $(Y, Z)$ by the solution $X^N$ of the forward equation

$$dX_t^N = \frac{1}{N} \sum_{j=1}^N \sigma(t, X_t^N, X_t^{N,j}) \, dW_t + \frac{1}{N} \sum_{j=1}^N b(t, X_t^N, X_t^{N,j}) \, dt,$$

(1.4)

$$X_0^N = x_0, \qquad t \in [0, T].$$

We point out the conceptional difference between our approximation and the standard approximation for the McKean–Vlasov SDE which considers as $X^N$ the dynamics $X'^{N,1}$ of the first element of a system of $N$ particles which dynamics is described by the system of SDEs

$$dX_t'^{N,i} = \frac{1}{N} \sum_{j=1}^N \sigma(t, X_t'^{N,i}, X_t'^{N,j}) \, dW_t^i + \frac{1}{N} \sum_{j=1}^N b(t, X_t'^{N,i}, X_t'^{N,j}) \, dt,$$

(1.5)

$$X_0^{N,i} = x_0, \qquad t \in [0, T], \ 1 \le i \le N,$$

governed by the independent Brownian motions $W^1, \ldots, W^N$ (see, e.g., Bossy [1], Bossy and Talay [2], Chan [4], Méléard [8], Overbeck [10], Pra and Hollander [14], Sznitman [15, 16], Talay and Vaillant [17]).

By proving that the speed of the convergence of $\Lambda^N = (X^N, Y^N, Z^N)$ to $\Lambda = (X, Y, Z)$ is of order $1/\sqrt{N}$ (see Proposition 3.2 and Theorem 3.2) we know that the sequence $\sqrt{N}(X^N - X, Y^N - Y, Z^N - Z)$, $N \ge 1$, is bounded,

$$\sup_{N \ge 1} E\left[ \sup_{t \in [0, T]} (|X_t^N - X_t|^2 + |Y_t^N - Y_t|^2) + \int_0^T |Z_t^N - Z_t|^2 \, dt \right] < +\infty,$$

so that the question of the limit behavior of the sequence $\sqrt{N}(X^N - X, Y^N - Y, Z^N - Z)$, $N \ge 1$, arises. Our special choice of the approximating sequence turns out to be particularly helpful for the study of this question. We show that the sequence $\sqrt{N}(X^N - X, Y^N - Y, Z^N - Z)$ converges in law on the space $C([0, T]; \mathbb{R}^d) \times C([0, T]; \mathbb{R}) \times L^2[0, T]; \mathbb{R}^d)$ to the unique solution $(\overline{X}, \overline{Y}, \overline{Z})$ of a (decoupled) forward–backward SDE of mean-field type

$$\overline{X}_t = \int_0^t \xi_s^{(1)}(X_s) \, ds + \int_0^t \xi_s^{(2)}(X_s) \, dW_s$$

(1.6)

$$+ \int_0^t (E[(\nabla_x b)(x, X_s)]_{x=X_s} \overline{X}_s + E[(\nabla_{x'} b)(x, X_s) \overline{X}_s]_{x=X_s}) \, ds$$

$$+ \int_0^t (E[(\nabla_x \sigma)(x, X_s)]_{x=X_s} \overline{X}_s + E[(\nabla_{x'} \sigma)(x, X_s) \overline{X}_s]_{x=X_s}) \, dW_s,$$

$$\overline{Y}_t = \{\xi^{(3)}(X_T) + E[\nabla_x \Phi(x, X_T)]_{x=X_T} \overline{X}_T + E[\nabla_{x'} \Phi(x, X_T) \overline{X}_T]_{x=X_T}\}$$

$$+ \int_t^T (\xi_s^{(4)}(\Lambda_s) + E[\nabla_\lambda f(\lambda, \Lambda_s)]_{\lambda=\Lambda_s} \overline{\Lambda}_s)$$

(1.7)
$$+ E[\nabla_{\lambda'} f(\lambda, \Lambda_s)\overline{\Lambda_s}]_{\lambda=\Lambda_s}) \, ds$$
$$- \int_t^T \overline{Z}_s \, dW_s,$$

which is not only driven by a $d$-dimensional Brownian motion $W$ but also by a zero-mean $(2d+1)$-parameter Gaussian field $\xi = \{\xi_t^{(1)}(x), \xi_t^{(2)}(x), \xi^{(3)}(x), \xi_t^{(4)}(x, y, z), (t, x, y, z) \in [0, T] \times \mathbb{R}^d \times \mathbb{R} \times \mathbb{R}^d\}$ which is independent of $W$ and whose covariance function is defined in Theorem 4.2.

Such a characterization, and in particular the study of convergence in law of a sequence of solutions of BSDEs to the solution of another BSDE, which includes not only the $Y$-component but also the $Z$-component and is not embedded in a Markovian framework is new to our best knowledge (observe that the solutions of SDEs of mean-field type are not Markovian). The main idea for the proof of the convergence result consists in rewriting the forward–backward equations for $\Lambda^N = (X^N, Y^N, Z^N)$ by using the stochastic field $\xi^N = (\xi^{1,N}, \xi^{2,N}, \xi^{3,N}, \xi^{4,N})$,

$$\xi_t^{i,N}(x) = \frac{1}{\sqrt{N}} \sum_{k=1}^N (\gamma(x, X_t^{N,k}) - E[\gamma(x, X_t^{N,k})]), \qquad (t, x) \in [0, T] \times \mathbb{R}^d$$

with $\gamma = b$ for $i = 1$, and $\gamma = \sigma$ for $i = 2$,

$$\xi^{3,N}(x) = \frac{1}{\sqrt{N}} \sum_{k=1}^N (\Phi(x, X_T^{N,k}) - E[\Phi(x, X_T^{N,k})]), \qquad x \in \mathbb{R}^d,$$

$$\xi_t^{4,N}(\lambda) = \frac{1}{\sqrt{N}} \sum_{k=1}^N (f(\lambda, \Lambda_t^{N,k}) - E[f(\lambda, \Lambda_t^{N,k})]), \qquad (t, \lambda) \in [0, T] \times \mathbb{R}^{2d+1},$$

and in proving that this stochastic field converges in law to $\xi$. The proof of this convergence in law on $C([0, T] \times \mathbb{R}^{2d+1}; \mathbb{R}^{d+d^2+2})$ (Proposition 4.3) is split into the proof of the tightness and that of the convergence of the finite-dimensional laws. For getting the tightness of the laws of $\xi^N, N \geq 1$, we have to suppose that the function $f((x, y, z), (x', y', z'))$ does not depend on $z'$, that is, there is no averaging with respect to the $Z$-component of the BSDE; in the proof of the convergence of the finite-dimensional laws the central limit theorem plays a crucial role. Once having the convergence in law of $\xi^N$ to $\xi$ we can apply Skorohod's representation theorem in order to obtain the almost sure convergence $\xi'^N \to \xi'$ for copies of $\xi^N$ and $\xi$ on an appropriate probability space. This allows to show that the associated redefinitions of the triplets $\sqrt{N}(X^N - X, Y^N - Y, Z^N - Z), N \geq 1$, converge on this new probability space to the solution $(\widehat{X}, \widehat{Y}, \widehat{Z})$ of the redefinition of the system (1.6)–(1.7) (Proposition 4.4). Our main result follows then easily.

Our paper is organized as follows: the short Section 2 recalls briefly some elements of the theory of backward SDEs which will be needed in what follows. In Section 3 we introduce the notion of MFBSDEs and the framework in which we study them, and we prove the existence and uniqueness. Moreover, the approximation of MFBSDE (1.1) by (1.2) is studied, and in a "Markovian-like" framework in which (1.1) and (1.2) are associated with the forward (1.3) and (1.4), respectively, the convergence speed is estimated. Finally, Section 4 is devoted to the study of the limit behavior of the triplet $\sqrt{N}(X^N - X, Y^N - Y, Z^N - Z)$, $N \geq 1$. In order to give to our approach a better readability we discuss first the limit behavior of the triplet $\sqrt{N}(X^N - X)$, $N \geq 1$, to investigate later the limit behavior of the triplets by using analogies of the arguments.

**2. Preliminaries.** The purpose of this section is to introduce some basic notions and results concerning BSDEs, which will be needed in the subsequent sections. In all that follows we will work on a slight extension of the classical Wiener space $(\Omega, \mathcal{F}, P)$:

• For an arbitrarily fixed time horizon $T > 0$ and a countable index set $I$ (which will be specified later), $\Omega$ is the set of all families $(\omega^i)_{i \in I}$ of continuous functions $\omega^i : [0, T] \to \mathbb{R}^d$ with initial value 0 $(\Omega = C_0([0, T]; \mathbb{R}^d)^I)$; it is endowed with the product topology generated by the uniform convergence on its components $C_0([0, T]; \mathbb{R}^d)$;

• $\mathcal{B}(\Omega)$ denotes the Borel $\sigma$-field over $\Omega$ and $B = (W^i)_{i \in I}$ is the coordinate process over $\Omega$: $W_t^i(\omega) = \omega_t^i, t \in [0, T], \omega \in \Omega, i \in I$;

• $P$ is the Wiener measure over $(\Omega, \mathcal{B}(\Omega))$, that is, the unique probability measure with respect to which $W^i, i \in I$, form a family of independent $d$-dimensional Brownian motions. Finally,

• $\mathcal{F}$ is the $\sigma$-field $\mathcal{B}(\Omega)$ completed with respect to the Wiener measure $P$. Let $W := W^0$. We endow our probability space $(\Omega, \mathcal{F}, P)$ with the filtration $\mathbf{F} = (\mathcal{F}_t)_{t \in [0, T]}$ which is generated by the Brownian motion $W$, enlarged by the $\sigma$-field $\mathcal{G} = \sigma\{W_t^i, t \in [0, T], i \in I \setminus \{0\}\}$ and completed by the collection $\mathcal{N}_P$ of all $P$-null sets:

$$\mathcal{F}_t = \mathcal{F}_t^W \vee \mathcal{G}, \qquad t \in [0, T],$$

where $\mathbf{F}^W = (\mathcal{F}_t^W = \sigma\{W_r, r \leq t\} \vee \mathcal{N}_P)_{t \in [0, T]}$. We observe that the Brownian motion $W$ has with respect to the filtration $\mathbf{F}$ the martingale representation property, that is, for every $\mathcal{F}_T$-measurable, square integrable random variable $\xi$ there is some $d$-dimensional $\mathbf{F}$-progressively measurable, square integrable process $Z = (Z_t)_{t \in [0, T]}$ such that

$$\xi = E[\xi | \mathcal{G}] + \int_0^T Z_t \, dW_t, \qquad P\text{-a.s.}$$

We also shall introduce the following spaces of processes which will be used frequently in the sequel:

$$S_{\mathbf{F}}^2([0,T]) = \Big\{ (Y_t)_{t \in [0,T]} \text{ continuous adapted process:}$$

$$E\Big[ \sup_{t \in [0,T]} |Y_t|^2 \Big] < +\infty \Big\};$$

$$L_{\mathbf{F}}^2([0,T]; \mathbb{R}^d) = \Big\{ (Z_t)_{t \in [0,T]} \mathbb{R}^d\text{-valued progressively measurable process:}$$

$$E\Big[ \int_0^T |Z_t|^2 \, dt \Big] < +\infty \Big\}.$$

(Recall that $|z|$ denotes the Euclidean norm of $z \in \mathbb{R}^n$.) Let us now consider a measurable function $g: \Omega \times [0,T] \times \mathbb{R} \times \mathbb{R}^d \to \mathbb{R}$ with the property that $(g(t,y,z))_{t \in [0,T]}$ is $\mathbf{F}$-progressively measurable for all $(y,z)$ in $\mathbb{R} \times \mathbb{R}^d$. We make the following standard assumptions on the coefficient $g$:

(A1) *There is some real $C \geq 0$ such that, $P$-a.s., for all $t \in [0,T], y_1, y_2 \in \mathbb{R}, z_1, z_2 \in \mathbb{R}^d$,*

$$|g(t,y_1,z_1) - g(t,y_2,z_2)| \leq C(|y_1 - y_2| + |z_1 - z_2|).$$

(A2) $g(\cdot, 0, 0) \in L_{\mathbf{F}}^2([0,T]; \mathbb{R})$.

The following result on BSDEs is by now well known; for its proof the reader is referred, for instance, to the pioneering work by Pardoux and Peng [11], but also to El Karoui, Peng and Quenez [5].

LEMMA 2.1. *Let the coefficient $g$ satisfy the assumptions* (A1) *and* (A2). *Then, for any random variable $\xi \in L^2(\Omega, \mathcal{F}_T, P)$, the BSDE associated with the data couple $(g, \xi)$*

$$Y_t = \xi + \int_t^T g(s, Y_s, Z_s) \, ds - \int_t^T Z_s \, dB_s, \qquad 0 \leq t \leq T,$$

*has a unique $\mathbf{F}$-progressively measurable solution*

$$(Y, Z) \in \mathcal{S}_{\mathbf{F}}^2([0,T]) \times L_{\mathbf{F}}^2([0,T]; \mathbb{R}^d).$$

The proof of Lemma 2.1 is related with the following, by now standard estimate for BSDEs.

LEMMA 2.2. *Let $(g_1, \xi_1), (g_2, \xi_2)$ be two data couples for which we suppose that $g_k$ satisfies the assumptions* (A.1) *and* (A.2) *and $\xi_k \in L^2(\Omega, \mathcal{F}_T, P)$, $k = 1, 2$. We denote by $(Y^k, Z^k)$ the unique solution of the BSDE with the data $(g_k, \xi_k), k = 1, 2$. Then, for every $\delta > 0$, there exist some $\gamma(= \gamma_\delta) > 0$* *and some $C(=C_\delta) > 0$ only depending on $\delta$ and on the Lipschitz constants of $g_k, k = 1, 2$, such that, with the notation*

$$(\overline{Y}, \overline{Z}) = (Y^1 - Y^2, Z^1 - Z^2), \qquad \overline{g} = g_1 - g_2, \qquad \overline{\xi} = \xi_1 - \xi_2,$$

*we have*

$$E\left[\int_0^T e^{\gamma t}(|\overline{Y}_t|^2 + |\overline{Z}_t|^2)\, dt\right] \leq CE[e^{\gamma T}|\overline{\xi}|^2] + \delta E\left[\int_0^T e^{\gamma t}|\overline{g}(t, Y_t^1, Z_t^1)|^2\, dt\right].$$

Besides the existence and uniqueness result we shall also recall the comparison theorem for BSDEs (see Theorem 2.2 in El Karoui, Peng and Quenez [5] or also Proposition 2.4 in Peng [13]).

LEMMA 2.3 (Comparison theorem). *Given two coefficients $g_1$ and $g_2$ satisfying* (A1) *and* (A2) *and two terminal values $\xi_1, \xi_2 \in L^2(\Omega, \mathcal{F}_T, P)$, we denote by $(Y^1, Z^1)$ and $(Y^2, Z^2)$ the solution of the BSDE with the data $(\xi_1, g_1)$ and $(\xi_2, g_2)$, respectively. Then we have:*

(i) (Monotonicity). *If $\xi_1 \geq \xi_2$ and $g_1 \geq g_2, a.s.$, then $Y_t^1 \geq Y_t^2$, for all $t \in [0, T]$, a.s.*

(ii) (Strict monotonicity). *If, in addition to* (i), *also $P\{\xi_1 > \xi_2\} > 0$, then we have $P\{Y_t^1 > Y_t^2\} > 0$, for all $0 \leq t \leq T$, and in particular, $Y_0^1 > Y_0^2$.*

After this short and very basic recall on BSDEs let us now investigate the limit approach for mean-field BSDEs (MFBSDEs).

## 3. Mean-field BSDEs.

3.1. *The notion of mean-field BSDEs. Existence and uniqueness.* The objective of our paper is to discuss a special mean-field problem in a purely stochastic approach. Let us first introduce the framework in which we want to study the limit approach for MFBSDEs. For this we specify the countable index set introduced in the proceeding section as follows:

$$I := \{i \mid i \in \{1, 2, 3, \ldots\}^k, k \geq 1\} \cup \{0\}.$$

For two elements $i = (i_1, \ldots, i_k)$, $i' = (i_1', \ldots, i_{k'}')$ of $I$ we define $i \oplus i' = (i_1, \ldots, i_k, i_1', \ldots, i_{k'}') \in I$, [with the convention that $(0) \oplus i = i$]. Then, in particular, for all $\ell \geq 1$, $(\ell) \oplus i = (\ell, i_1, \ldots, i_k)$.

We also shall introduce a family of shift operators $\Theta^k : \Omega \to \Omega$, $k \geq 0$, over $\Omega$. For this end we set $\Theta^k(\omega) = (\omega^{(k) \oplus j})_{j \in I}, \omega \in \Omega, k \geq 0$, and we observe that $\Theta^k(\omega)$ can be regarded as an element of $\Omega$ and $\Theta^k$ as an operator mapping $\Omega$ into $\Omega$. The fact that all these operators $\Theta^k : \Omega \to \Omega$ let the Wiener measure $P$ invariant (i.e., $P_{\Theta^k} = P$) allows to interpret $\Theta^k$ as operator defined over $L^0(\Omega, \mathcal{F}, P)$: putting $\Theta^k(\xi)(\omega) := \xi(\Theta^k(\omega)), \omega \in \Omega$, for the random variables of the form $\xi(\omega) = f(\omega_{t_1}^{i_1}, \ldots, \omega_{t_n}^{i_n})$, $i_1, \ldots, i_n \in I, t_1, \ldots, t_n \in [0, T], f \in C(R^{d \times n}), n \geq 1$, we can extend this definition from this class of continuous Wiener functionals to the whole space $L^0(\Omega, \mathcal{F}, P)$ by using the density of the class of smooth Wiener functionals in $L^0(\Omega, \mathcal{F}, P)$. We observe that, for all $\xi \in L^0(\Omega, \mathcal{F}, P)$, the random variables $\Theta^k(\xi), k \geq 1$, are independent and identically distributed (i.i.d.), of the same law as $\xi$ and independent of the Brownian motion $W$.

Finally, to shorten notation we introduce the $(N + 1)$-dimensional shift operator $\Theta_N = (\Theta^0, \Theta^1, \ldots, \Theta^N)$, which associates a random variable $\xi \in L^0(\Omega, \mathcal{F}, P)$ with the $(N+1)$-dimensional random vector $\Theta_N(\xi) = (\xi, \Theta^1(\xi), \ldots, \Theta^N(\xi))$ (notice that $\Theta^0$ is the identical operator). If $\xi$ on its part is a random vector, $\Theta^k(\xi)$ and $\Theta_N(\xi)$ are defined by a componentwise application of the corresponding operators.

For an arbitrarily fixed natural number $N \geq 0$ let us now consider a measurable function $f : \Omega \times [0, T] \times \mathbb{R}^{N+1} \times \mathbb{R}^{(N+1) \times d} \to \mathbb{R}$ with the property that $(f(t, \mathbf{y}, \mathbf{z}))_{t \in [0,T]}$ is $\mathbf{F}$-progressively measurable for all $(\mathbf{y}, \mathbf{z})$ in $\mathbb{R}^{N+1} \times \mathbb{R}^{(N+1) \times d}$. We make the following standard assumptions on the coefficient $f$, which extend (A1) and (A2) in a natural way:

(B1) *There is some constant $C \geq 0$ such that, $P$-a.s., for all $t \in [0, T], \mathbf{y}_1,$*
*$\mathbf{y}_2 \in \mathbb{R}^{N+1}, \mathbf{z}_1, \mathbf{z}_2 \in \mathbb{R}^{(N+1) \times d}$,*

$$|f(t, \mathbf{y}_1, \mathbf{z}_1) - f(t, \mathbf{y}_2, \mathbf{z}_2)| \leq C(|\mathbf{y}_1 - \mathbf{y}_2| + |\mathbf{z}_1 - \mathbf{z}_2|).$$

(B2) $f(\cdot, 0, 0) \in L_{\mathbf{F}}^2(0, T; \mathbb{R})$.

The following proposition extends Lemma 2.1 to the type of backward equations which will be used for the approximation of MFBSDEs.

PROPOSITION 3.1. *Let the function $f$ satisfy the above assumptions (B1) and (B2). Then, for any random variable $\xi \in L^2(\Omega, \mathcal{F}_T, P)$, the BSDE associated with $(f, \xi)$*

$$(3.1) \qquad dY_t = -f(t, \Theta_N(Y_t, Z_t)) \, dt + Z_t \, dW_t, \qquad t \in [0, T], \qquad Y_T = \xi$$

*has a unique adapted solution $(Y, Z) \in S_{\mathbf{F}}^2([0, T]; \mathbb{R}) \times L_{\mathbf{F}}^2([0, T]; \mathbb{R}^d)$.*

PROOF. Let $H^2 := L_{\mathbf{F}}^2([0, T]; \mathbb{R}) \times L_{\mathbf{F}}^2([0, T]; \mathbb{R}^d)$. It is sufficient to prove the existence and the uniqueness for the above BSDE in $H^2$. Indeed, if $(Y, Z)$ is a solution of our BSDE in $H^2$ an easy standard argument shows that it is also in $B^2 := S_{\mathbf{F}}^2([0, T]; \mathbb{R}) \times L_{\mathbf{F}}^2([0, T]; \mathbb{R}^d)$. On the other hand, the uniqueness in $H^2$ implies obviously that its subspace $B^2$.

For proving the existence and uniqueness in $H^2$ we consider for an arbitrarily given couple of processes $(U, V) \in H^2$ the coefficient $g_t = f(t, \Theta_N(U_t,$ $V_t$)), $t \in [0, T]$. Since $g$ is an element of $L^2_{\mathbf{F}}([0, T]; \mathbb{R})$ it follows from Lemma 2.1 that there is a unique solution $\Phi(U, V) := (Y, Z) \in H^2$ of the BSDE

$$dY_t = -f(t, \Theta_N(U_t, V_t)) \, dt + Z_t \, dW_t, \qquad t \in [0, T], \qquad Y_T = \xi.$$

For a such defined mapping $\Phi \colon H^2 \to H^2$ it suffices to prove that it is a contraction with respect to an appropriate equivalent norm on $H^2$ in order to complete the proof. For this end we consider two couples $(U^1, V^1), (U^2, V^2) \in H$ and $(Y^k, Z^k) = \Phi(U^k, V^k), k = 1, 2$. Then, due to Lemma 2.2, for all $\delta > 0$ there is some constant $\gamma > 0$ (only depending on $\delta$) such that, with the notation $(\overline{Y}, \overline{Z}) = (Y^1 - Y^2, Z^1 - Z^2)$,

$$E\left[\int_0^T e^{\gamma t}(|\overline{Y}_t|^2 + |\overline{Z}_t|^2) \, dt\right]$$
$$\leq \delta E\left[\int_0^T e^{\gamma t} |f(t, \Theta_N(U_t^1, V_t^1)) - f(t, \Theta_N(U_t^2, V_t^2))|^2 \, dt\right].$$

Let $(\overline{U}, \overline{V}) = (U^1 - U^2, V^1 - V^2)$. Then, from the Lipschitz continuity of $f(\omega, t, \cdot, \cdot)$ [with Lipschitz constant $L$ which does not depend on $(\omega, t)$] and the fact that the random vectors $\Theta^k(\overline{U}, \overline{V}) (= \Theta^k(U^1, V^1) - \Theta^k(U^2, V^2)), 0 \leq k \leq N$ obey all the same law, we have

$$E\left[\int_0^T e^{\gamma t}(|\overline{Y}_t|^2 + |\overline{Z}_t|^2) \, dt\right] \leq \delta L^2 E\left[\int_0^T e^{\gamma t}\left|\sum_{k=0}^N |\Theta^k(\overline{U}_t, \overline{V}_t)|\right|^2 \, dt\right]$$
$$\leq \delta L^2 (N+1) \int_0^T e^{\gamma t} \sum_{k=0}^N E[|\Theta^k(\overline{U}_t, \overline{V}_t)|^2] \, dt$$
$$= \delta L^2 (N+1)^2 \int_0^T e^{\gamma t} E[|\overline{U}_t|^2 + |\overline{V}_t|^2] \, dt$$
$$= \frac{1}{2} E\left[\int_0^T e^{\gamma t}(|\overline{U}_t|^2 + |\overline{V}_t|^2) \, dt\right]$$

for $\delta := \frac{1}{2} L^{-2} (N+1)^{-2}$. This shows that if we endow the space $H^2$ with the norm

$$\|(U, V)\|_{H^2} = \left(E\left[\int_0^T e^{\gamma t}(|U_t|^2 + |V_t|^2) \, dt\right]\right)^{1/2}, \qquad (U, V) \in H^2,$$

the mapping $\Phi \colon H^2 \to H^2$ becomes a contraction. Thus, the proof is complete. $\square$

We now introduce the framework for the study of the limit of the above BSDE as $N$ tends to $+\infty$. For this end let be given a data triplet $(\Phi, g, \mathcal{X})$ with the following properties (C1), (C2) and (C3):

(C1) $g\colon \Omega \times [0,T] \times (\mathbb{R}^m \times \mathbb{R} \times \mathbb{R}^d)^2 \to \mathbb{R}$ is a bounded measurable function which is Lipschitz in $(\mathbf{u}, \mathbf{v}) \in (\mathbb{R}^m \times \mathbb{R} \times \mathbb{R}^d)^2$ with a Lipschitz constant $C$, that is, $P$-a.s., for all $t \in [0,T]$ and $(\mathbf{u}, \mathbf{v}), (\mathbf{u}', \mathbf{v}') \in (\mathbb{R}^m \times \mathbb{R} \times \mathbb{R}^d)^2$,

$$|g(t,(\mathbf{u},\mathbf{v})) - g(t,(\mathbf{u}',\mathbf{v}'))| \le C(|\mathbf{u} - \mathbf{u}'| + |\mathbf{v} - \mathbf{v}'|);$$

(C2) $\Phi\colon \Omega \times \mathbb{R}^m \times \mathbb{R}^m \to \mathbb{R}$ is a bounded measurable function such that $\Phi(\omega, \cdot, \cdot)$ is Lipschitz with a Lipschitz constant $C$, that is, $P$-a.s., for all $(x, \hat{x})$, $(x', \hat{x}') \in \mathbb{R}^m$,

$$|\Phi(x,\hat{x}) - \Phi(x',\hat{x}')| \le C(|x - x'| + |\hat{x} - \hat{x}'|).$$

(C3) $\mathcal{X} = (X^N)_{N \ge 1}$ is a Cauchy sequence in $\mathcal{S}_{\mathbf{F}}^2([0,T]; \mathbb{R}^m)$, that is, there is a (unique) process $X \in \mathcal{S}_{\mathbf{F}}^2([0,T]; \mathbb{R}^m)$ such that

$$E\left[ \sup_{t \in [0,T]} |X_t^N - X_t|^2 \right] \to 0 \qquad \text{as } N \to +\infty.$$

REMARK 3.1. In the next section we will consider as sequence $X^N$, $N \ge 1$, a special approximation of the solution of a forward SDE of McKean–Vlasov type. This special choice will allow to study the convergence speed of the BSDEs as $N$ tends to $+\infty$ and to characterize the nature of this convergence more precisely.

Given such a triplet $(\Phi, g, \mathcal{X})$ satisfying the assumptions (C1)–(C3) we put

$$f^N(\omega, t, \mathbf{y}, \mathbf{z}) := \frac{1}{N} \sum_{k=1}^N g(\Theta^k(\omega), t, X_t^N(\omega), (y_0, z_0), X_t^N(\Theta^k(\omega)), (y_k, z_k))$$

for $(\omega, t) \in \Omega \times [0,T], \mathbf{y} = (y_0, \ldots, y_N) \in \mathbb{R}^{N+1}, \mathbf{z} = (z_0, \ldots, z_N) \in \mathbb{R}^{(N+1) \times d}$, and

$$\xi^N(\omega) := \frac{1}{N} \sum_{k=1}^N \Phi(\Theta^k(\omega), X_T^N(\omega), X_T^N(\Theta^k(\omega))), \qquad N \ge 1.$$

We observe that, for every $N \ge 1$, $f^N$ satisfies the assumptions (B1) and (B2). Thus, due to Proposition 3.1, for all $N \ge 1$, there is a unique solution $(Y^N, Z^N)$ of the BSDE($N$)

$$(3.2) \quad Y_t^N = \xi^N + \int_t^T f^N(s, \Theta_N(Y_s^N, Z_s^N))\, ds - \int_t^T Z_s^N\, dW_s, \qquad t \in [0,T].$$

We remark in particular that the driving coefficient of the above BSDE($N$) (3.2) can be written as follows:

$$f^N(s, \Theta_N(Y_s^N, Z_s^N)) = \frac{1}{N} \sum_{k=1}^N (\Theta^k g)(s, (X_s^N, Y_s^N, Z_s^N), \Theta^k(X_s^N, Y_s^N, Z_s^N)),$$

$s \in [0, T]$ [recall that $\Theta^k g(\omega, s, (\mathbf{u}, \mathbf{v})) := g(\Theta^k(\omega), s, (\mathbf{u}, \mathbf{v}))$]. Our objective is to show that the unique solution of BSDE($N$) (3.2) converges in $\mathcal{B}^2 = \mathcal{S}_{\mathbf{F}}^2([0, T]) \times L_{\mathbf{F}}^2([0, T]; \mathbb{R}^d)$ to the unique solution $(Y, Z)$ of the MFBSDE

(3.3)
$$dY_t = -E[g(t, \mathbf{u}, \Lambda_t)]_{|\mathbf{u} = \Lambda_t} \, dt + Z_t \, dW_t, \qquad t \in [0, T],$$
$$Y_T = E[\Phi(x, X_T)]_{|x = X_T},$$

where we have used the notation $\Lambda = (X, Y, Z)$.

LEMMA 3.1. *Under the assumptions* (C1)–(C3) *the above MFBSDE (3.3) possesses a unique solution* $(Y, Z) \in \mathcal{B}^2$.

Since the proof is straightforward and uses essentially the argument developed in the proof of Proposition 3.1, we omit it. See also [3].

We now can formulate the following theorem:

THEOREM 3.1. *Under the assumptions* (C1)–(C3) *the unique solution* $(Y^N, Z^N)$ *of BSDE(N) (3.2) converges in* $\mathcal{B}^2$ *to the unique solution* $(Y, Z)$ *of the above MFBSDE (3.3):*

$$E\left[ \sup_{t \in [0, T]} |Y_t^N - Y_t|^2 + \int_0^T |Z_t^N - Z_t|^2 \, dt \right] \to 0 \qquad as \ N \to +\infty.$$

PROOF.  We first notice that

STEP 1.  For all $p \geq 2$,

$$E\left[ \int_0^T \left| \frac{1}{N} \sum_{k=1}^N (\Theta^k g)(t, \Lambda_t, \Theta^k(\Lambda_t)) - E[g(t, \mathbf{u}, \Lambda_t)]_{|\mathbf{u} = \Lambda_t} \right|^p dt \right] \to 0;$$

$$E\left[ \left| \frac{1}{N} \sum_{k=1}^N (\Theta^k \Phi)(X_T, \Theta^k(X_T)) - E[\Phi(x, X_T)]_{|x = X_T} \right|^p \right] \to 0$$

as $N \to +\infty$ [recall that $(\Theta^k \Phi)(\omega, X_T, \Theta^k(X_T))(\omega) := \Phi(\Theta^k(\omega), X_T(\omega), X_T(\Theta^k(\omega)))$].

For proving the first convergence we consider arbitrary $t \in [0, T]$ and $\mathbf{u} \in \mathbb{R}^m \times \mathbb{R} \times \mathbb{R}^d$. Observing that the sequence of random variables $(\Theta^k g)(t, \mathbf{u}, \Lambda_t)$, $k \geq 1$, is i.i.d. and of the same law as $g(t, \mathbf{u}, \Lambda_t)$, we obtain from the Strong Law of Large Numbers that

$$\frac{1}{N} \sum_{k=1}^N (\Theta^k g)(t, \mathbf{u}, \Lambda_t) \longrightarrow E[g(t, \mathbf{u}, \Lambda_t)],$$

$P$-a.s., as $N \to +\infty$. Let now, for an arbitrarily small $\varepsilon > 0$, $\Lambda_t^\varepsilon : \Omega \to \mathbb{R}^m \times \mathbb{R} \times \mathbb{R}^d$ be a random vector which has only a countable number of values and is such that $|\Lambda_t - \Lambda_t^\varepsilon| \le \varepsilon$, everywhere on $\Omega$. Then, obviously,

$$\frac{1}{N} \sum_{k=1}^N (\Theta^k g)(t, \Lambda_t^\varepsilon, \Theta^k(\Lambda_t)) \longrightarrow E[g(t, \mathbf{u}, \Lambda_t)]_{|\mathbf{u} = \Lambda_t^\varepsilon},$$

$P$-a.s., as $N$ tends to $+\infty$, and from the Lipschitz continuity of $g(\omega, t, \cdot, \mathbf{v})$, which is uniform in $(\omega, t, \mathbf{v})$, it follows that we also have the convergence for $\Lambda_t$ instead of $\Lambda_t^\varepsilon$:

$$\frac{1}{N} \sum_{k=1}^N (\Theta^k g)(t, \Lambda_t, \Theta^k(\Lambda_t)) \longrightarrow E[g(t, \mathbf{u}, \Lambda_t)]_{|\mathbf{u} = \Lambda_t},$$

$P$-a.s., as $N \to +\infty$. Finally, in view of the boundedness of $g$ and, hence, of that of the convergence, we get the announced result. With the same argument we also get the $L^p$-convergence for the terminal condition, for all $p \ge 2$.

STEP 2. By applying Lemma 2.2 to estimate the distance between the solution $(Y^N, Z^N)$ of BSDE($N$) (3.2) and that of the MFBSDE (3.3) we get, with the notation $\Lambda = (X, Y, Z)$ and $\Lambda^N = (X^N, Y^N, Z^N)$, that for any $\delta \in (0, 1)$ there are some $\gamma > 0$ and $C > 0$ only depending on $\delta$ (and, hence, in particular independent of $N$) such that

$$E\left[ \int_0^T e^{\gamma t} (|Y_t^N - Y_t|^2 + |Z_t^N - Z_t|^2) \, dt \right]$$

$$\le CE\left[ e^{\gamma T} \left| \frac{1}{N} \sum_{k=1}^N (\Theta^k \Phi)(X_T^N, \Theta^k(X_T^N)) - E[\Phi(x, X_T)]_{|x = X_T} \right|^2 \right]$$

$$+ \delta E\left[ \int_0^T e^{\gamma t} \left| \frac{1}{N} \sum_{k=1}^N (\Theta^k g)(t, \Lambda_t^N, \Theta^k(\Lambda_t^N)) \right.\right.$$

$$\left.\left. - E[g(t, \mathbf{u}, \Lambda_t)]_{|\mathbf{u} = \Lambda_t} \right|^2 dt \right].$$

Hence, taking into account that $g(\omega, t, \cdot, \cdot)$ and $\Phi(\omega, \cdot, \cdot)$ are Lipschitz with a Lipschitz constant $L > 0$ which does not depend on $(\omega, t)$, we deduce that

$$E\left[ \int_0^T e^{\gamma t} (|Y_t^N - Y_t|^2 + |Z_t^N - Z_t|^2) \, dt \right]$$

$$\le R_N + 8CL^2 E[e^{\gamma T} |X_T^N - X_T|^2]$$

$$+ 8L^2 \delta E\left[ \int_0^T e^{\gamma t} (|X_t^N - X_t|^2 + |Y_t^N - Y_t|^2 + |Z_t^N - Z_t|^2) \, dt \right.$$

with

$$R_N = 2CE\left[e^{\gamma T}\left|\frac{1}{N}\sum_{k=1}^{N}(\Theta^k\Phi)(X_T,\Theta^k(X_T)) - E[\Phi(x,X_T)]_{|x=X_T}\right|^2\right]$$

$$+ 2\delta E\left[\int_0^T e^{\gamma t}\left|\frac{1}{N}\sum_{k=1}^{N}(\Theta^k g)(t,\Lambda_t,\Theta^k(\Lambda_t)) - E[g(t,\mathbf{u},\Lambda_t)]_{|\mathbf{u}=\Lambda_t}\right|^2 dt\right].$$

Consequently, for $\delta := (16L^2)^{-1}$,

$$\frac{1}{2}E\left[\int_0^T e^{\gamma t}(|Y_t^N - Y_t|^2 + |Z_t^N - Z_t|^2)\,dt\right]$$

$$\leq R_N + 8CL^2E[e^{\gamma T}|X_T^N - X_T|^2] + \frac{1}{2}E\left[\int_0^T e^{\gamma t}|X_t^N - X_t|^2\,dt\right].$$

From Step 1 and (C3) it then follows that

$$E\left[\int_0^T e^{\gamma t}(|Y_t^N - Y_t|^2 + |Z_t^N - Z_t|^2)\,dt\right] \to 0 \qquad \text{as } N \to +\infty.$$

For proving that $(Y^N, Z^N)$ converges to $(Y,Z)$ also in $\mathcal{B}^2$ we have still to show the convergence of $Y^N$ to $Y$ in $\mathcal{S}_{\mathbf{F}}^2([0,T])$. To this end we apply Itô's formula to $|Y_t^N - Y_t|^2$ and take then the conditional expectation. Thus, because of the boundedness of $g$, we have:

$$|Y_t^N - Y_t|^2 + E\left[\int_t^T |Z_s^N - Z_s|^2\,ds|\mathcal{F}_t\right]$$

$$= E\left[\left|\frac{1}{N}\sum_{k=1}^{N}(\Theta^k\Phi)(X_T^N,\Theta^k(X_T^N)) - E[\Phi(x,X_T)]_{|x=X_T}\right|^2\Big|\mathcal{F}_t\right]$$

$$+ 2E\left[\int_t^T (Y_s^N - Y_s)\right.$$

$$\left. \times\left\{\frac{1}{N}\sum_{k=1}^{N}(\Theta^k g)(s,\Lambda_s^N,\Theta^k(\Lambda_s^N)) - E[g(s,\mathbf{u},\Lambda_s)]_{|\mathbf{u}=\Lambda_s}\right\}ds|\mathcal{F}_t\right]$$

$$\leq E\left[\left|\frac{1}{N}\sum_{k=1}^{N}(\Theta^k\Phi)(X_T^N,\Theta^k(X_T^N)) - E[\Phi(x,X_T)]_{|x=X_T}\right|^2\Big|\mathcal{F}_t\right]$$

$$+ CE\left[\int_0^T |Y_s^N - Y_s|\,ds|\mathcal{F}_t\right], \qquad t \in [0,T].$$

Therefore,

$$E\left[\sup_{t\in[0,T]} |Y_t^N - Y_t|^4\right]$$

$$\leq 4E\left[\left|\frac{1}{N}\sum_{k=1}^{N}(\Theta^k\Phi)(X_T^N,\Theta^k(X_T^N))-E[\Phi(x,X_T)]_{|x=X_T}\right|^4\right]$$

$$+4C^2TE\left[\int_0^T|Y_s^N-Y_s|^2\,ds\right]$$

and the announced convergence follows now from the preceding result and from Step 1. $\square$

3.2. *Convergence speed of the approximation of the MFBSDE.* In Theorem 3.1 we have seen that under the assumptions (C1)–(C3) the solution $(Y^N,Z^N)$ of BSDE($N$) (3.2) converges toward the unique solution $(Y,Z)$ of our MFBSDE (3.3). The objective of this section is to study the speed of this convergence in the special case where the sequence $(X^N)_{N\geq 1}$ is an approximation of a forward SDE of McKean–Vlasov type. For this let $b\colon\Omega\times[0,T]\times\mathbb{R}^m\times\mathbb{R}^m\to\mathbb{R}^m$ and $\sigma\colon\Omega\times[0,T]\times\mathbb{R}^m\times\mathbb{R}^m\to\mathbb{R}^{m\times d}$ be bounded measurable functions which are supposed to be Lipschitz in $(x,x')$ with a Lipschitz constant $C$, that is, $P$-a.s.,

(D) $|b(t,x,x')-b(t,\hat{x},\hat{x}')|+|\sigma(t,x,x')-\sigma(t,\hat{x},\hat{x}')|\leq C(|x-\hat{x}|+|x'-\hat{x}'|)$, for all $t\in[0,T],(x,x'),(\hat{x},\hat{x}')\in\mathbb{R}^m\times\mathbb{R}^m$.

We consider the following forward equation of McKean–Vlasov type

$$(3.4)\qquad\begin{aligned}dX_t&=E[b(t,x,X_t)]_{|x=X_t}\,dt+E[\sigma(t,x,X_t)]_{|x=X_t}\,dW_t,\\X_0&=x\in\mathbb{R}^m,\qquad t\in[0,T],\end{aligned}$$

and we approximate this equation by forward equations in the same spirit as we have approximated our MFBSDE (3.3) by BSDE($N$) (3.2), $N\geq 1$. More precisely, we consider as approximating SDE($N$):

$$(3.5)\qquad\begin{aligned}dX_t^N&=\frac{1}{N}\sum_{k=1}^{N}(\Theta^k b)(t,X_t^N,\Theta^k(X_t^N))\,dt\\&\quad+\frac{1}{N}\sum_{k=1}^{N}(\Theta^k\sigma)(t,X_t^N,\Theta^k(X_t^N))\,dW_t,\\X_0^N&=x,\qquad t\in[0,T].\end{aligned}$$

For the above forward equations we can state the following result.

PROPOSITION 3.2. *Under the above standard assumptions on the coefficients $b$ and $\sigma$ we have:*

(i) *The forward SDE of McKean–Vlasov type possesses a unique solution $X\in\mathcal{S}_{\mathbb{F}}^2([0,T];\mathbb{R}^m)$. Moreover, $X$ is adapted with respect to the filtration $\mathbf{F}^W$ generated by the Brownian motion $W$.*

(ii) *For all $N \geq 1$, the forward equation SDE(N) (3.5) admits a unique solution $X^N \in \mathcal{S}_{\mathbf{F}}^2([0,T]; \mathbb{R}^m)$.*

(iii)

$$E\left[\sup_{t \in [0,T]} |X_t^N - X_t|^2\right] \leq \frac{C}{N} \qquad \text{for all } N \geq 1,$$

*where $C$ is a real constant which only depends on the bounds and the Lipschitz constants of the coefficients $b$ and $\sigma$.*

PROOF. By using the properties of the operators $\Theta^k, k \geq 1$, the results (i) and (ii) can be obtained by easy standard estimates for SDEs. For this we see that the coefficients $\frac{1}{N}\sum_{k=1}^{N}(\Theta^k b)(t,x,x')$ and $\frac{1}{N}\sum_{k=1}^{N}(\Theta^k \sigma)(t,x,x')$ of SDE(N) (3.5) are bounded, $\mathbf{F}$-progressively measurable for all $(x,x')$, and Lipschitz in $(x,x')$, uniformly with respect to $(\omega,t) \in \Omega \times [0,T]$. Moreover, the bound and the Lipschitz constant do not depend on $N$. The coefficients of the SDE (3.4) of McKean–Vlasov type, too, are bounded and Lipschitz, uniformly with respect to $t \in [0,T]$. From the fact that the coefficients $E[b(t,x,X_t)]$ and $E[\sigma(t,x,X_t)]$ are deterministic it follows easily that the unique solution of the McKean–Vlasov type equation is $\mathbf{F}^W$-progressively measurable.

For the proof of statement (iii) of the proposition we notice that, thanks to the fact that the processes $\Theta^k(\sigma(\cdot,x,X)) = ((\Theta^k \sigma)(t,x,\Theta^k(X_t)))_{t \in [0,T]}, k \geq 1$, are mutually independent and identically distributed, of the same law as the process $\sigma(\cdot,x,X) = (\sigma(t,x,X_t))_{t \in [0,T]}$, and independent of $W$ (and, hence, also of $X$), we get

$$E\left[\left|\frac{1}{N}\sum_{k=1}^{N}(\Theta^k \sigma)(t,X_t,\Theta^k(X_t)) - E[\sigma(t,x,X_t)]_{|x=X_t}\right|^2\right]$$

$$= \int_{\mathbb{R}^m} E\left[\left|\frac{1}{N}\sum_{k=1}^{N}(\Theta^k \sigma)(t,x,\Theta^k(X_t)) - E[\sigma(t,x,X_t)]\right|^2\right] P_{X_t}(dx)$$

$$= \sum_{1 \leq i \leq m, 1 \leq j \leq d} \int_{\mathbb{R}^m} \text{Var}\left(\frac{1}{N}\sum_{k=1}^{N}(\Theta^k \sigma_{i,j})(t,x,\Theta^k(X_t))\right) P_{X_t}(dx)$$

$$= \frac{1}{N}\sum_{1 \leq i \leq m, 1 \leq j \leq d} \int_{\mathbb{R}^m} \text{Var}(\sigma_{i,j}(t,x,X_t)) P_{X_t}(dx) \leq \frac{C}{N}$$

for some constant $C$ which does not depend on $N \geq 1$. Similarly, we have

$$E\left[\left|\frac{1}{N}\sum_{k=1}^{N}(\Theta^k b)(t,X_t,\Theta^k(X_t)) - E[b(t,x,X_t)]_{|x=X_t}\right|^2\right] \leq \frac{C}{N}.$$

Using the fact that $\frac{1}{N}\sum_{k=1}^{N}(\Theta^{k}\sigma)(t,\cdot,\cdot)$ and $\frac{1}{N}\sum_{k=1}^{N}(\Theta^{k}b)(t,\cdot,\cdot)$ are Lipschitz, uniformly with respect to $t \in [0,T]$, and we obtain now statement (iii) by an SDE standard estimate. $\quad\square$

The above estimate of the convergence speed can be extended from the forward equations to the associated BSDEs. Indeed, we have the following.

THEOREM 3.2.  *We assume* (C1)–(C2) *and take instead of* (C3) *the assumption that $X^{N} \in \mathcal{S}^{2}([0,T];\mathbb{R}^{m})$ is the unique solution of SDE(N) (3.5), $N \geq 1$. Let $(Y,Z)$ denote the unique solution of MFBSDE (3.3) driven by the forward SDE (3.4) of McKean–Vlasov type with solution $X$, and let $(Y^{N}, Z^{N})$ be the unique solution of BSDE(N) (3.2) driven by SDE(N) (3.5). Then, for some constant $C$ which depends only on the bounds and the Lipschitz constants of $b, \sigma$ and $f$,*

$$E\left[\sup_{t\in[0,T]}|Y_{t}^{N}-Y_{t}|^{2} + \int_{0}^{T}|Z_{t}^{N}-Z_{t}|^{2}\right] \leq \frac{C}{N} \qquad \text{for all } N \geq 1.$$

PROOF.  We first notice that, since due to Proposition 3.2 the process $X$ is $\mathbf{F}^{W}$-progressively measurable, the unique solution $(Y,Z)$ of MFBSDE (3.3) driven by the process $X$ must also be $\mathbf{F}^{W}$-progressively measurable. Consequently, the triplet $\Lambda = (X,Y,Z)$ is independent of all the processes $\Theta^{k}(f(\cdot,\mathbf{u},\Lambda)) = ((\Theta^{k}f)(t,\mathbf{u},\Theta^{k}(\Lambda_{t})))_{t\in[0,T]}$, $k \geq 1$, which allows to use the argument developed in the proof of Proposition 3.2 and to show that, for some constant $C$ which only depends on the bound of $f$,

$$E\left[\left|\frac{1}{N}\sum_{k=1}^{N}(\Theta^{k}f)(t,\Lambda_{t},\Theta^{k}(\Lambda_{t})) - E[f(t,\mathbf{u},\Lambda_{t})]_{|\mathbf{u}=\Lambda_{t}}\right|^{2}\right] \leq \frac{C}{N}$$

$$\text{for all } N \geq 1.$$

Analogously, we get, again for some constant $C$ which only depends on the bound of $\Phi$,

$$E\left[\left|\frac{1}{N}\sum_{k=1}^{N}(\Theta^{k}\Phi)(X_{T},\Theta^{k}(X_{T})) - E[\Phi(x,X_{T})]_{|x=X_{T}}\right|^{2}\right] \leq \frac{C}{N} \qquad \text{for all } N \geq 1.$$

These estimates together with Proposition 3.2 and BSDE standard estimates yield the wished speed of the convergence of $(Y^{N}, Z^{N})$ to $(Y,Z)$. The proof is complete. $\quad\square$

## 4. A central limit theorem for MFBSDEs.

We have seen in the preceding Section 3.1 that, if the sequence $X^{N}, N \geq 1$, is defined by the forward equations SDE($N$) (3.5), the speed of the convergence of $(X^N, Y^N, Z^N)$ to $(X, Y, Z)$ is of order $1/\sqrt{N}$, that is, for some real $C$,

$$E\left[\sup_{t \in [0,T]} (|X_t^N - X_t|^2 + |Y_t^N - Y_t|^2) + \int_0^T |Z_t^N - Z_t|^2 \, dt\right] \le \frac{C}{N} \qquad \text{for } N \ge 1.$$

Consequently, the sequence of processes $\sqrt{N}(X^N - X, Y^N - Y, Z^N - Z)$, $N \ge 1$, is bounded in $\mathcal{S}_{\mathbf{F}}^2([0,T]; \mathbb{R}^m) \times \mathcal{S}_{\mathbf{F}}^2([0,T]) \times \mathcal{L}_{\mathbf{F}}^2([0,T]; \mathbb{R}^d)$ and the question about the limit behavior of this sequence arises. The study of this question is the objective of this section. For this end we will begin to investigate in a first section the limit behavior of the sequence $(\sqrt{N}(X^N - X))_{N \ge 1} \subset \mathcal{S}_{\mathbf{F}}^2([0,T]; \mathbb{R}^m)$. This discussion which involves, for its most part, some recall from known facts prepares the study of the limit behavior of the triplet $\sqrt{N}(X^N - X, Y^N - Y, Z^N - Z)$, $N \ge 1$, which will be the object of the second section.

For the sake of simplicity of the notation but without restricting the generality of our method we suppose in what follows that the coefficients $b, \sigma, f$ and $\Phi$ do not depend on $(\omega, t)$. Moreover, we will assume that the functions are continuously differentiable.

### 4.1. Limit behavior of $\sqrt{N}(X^N - X)$.

Let $b : \mathbb{R}^d \times \mathbb{R}^d \to \mathbb{R}^d$ and $\sigma : \mathbb{R}^d \times \mathbb{R}^d \to \mathbb{R}^{d \times d}$ be two bounded, continuously differentiable functions with bounded first-order derivatives, and let $x_0$ be an arbitrarily fixed element of $\mathbb{R}^d$. As in the preceding section, but now under the additional assumption that the coefficients do not depend on $(\omega, t)$, we denote by $X^N \in \mathcal{S}_{\mathbf{F}}^2([0,T]; \mathbb{R}^d)$ the unique solution of SDE($N$)

$$
\begin{aligned}
(4.1) \qquad X_t^N = x_0 &+ \int_0^t \frac{1}{N} \sum_{k=1}^N b(X_s^N, \Theta^k(X_s^N)) \, ds \\
&+ \int_0^t \frac{1}{N} \sum_{k=1}^N \sigma(X_s^N, \Theta^k(X_s^N)) \, dW_s, \qquad t \in [0,T]
\end{aligned}
$$

and $X \in \mathcal{S}_{\mathbf{F}^W}^2([0,T]; \mathbb{R}^d)$ is the unique solution of McKean–Vlasov equation

$$
\begin{aligned}
(4.2) \qquad X_t = x_0 &+ \int_0^t E[b(x, X_s)]_{|x=X_s} \, ds + \int_0^t E[\sigma(x, X_s)]_{|x=X_s} \, dW_s, \\
&\qquad\qquad\qquad\qquad\qquad\qquad\qquad\qquad s \in [0,T].
\end{aligned}
$$

Then we have the limit behavior of $(\sqrt{N}(X^N - X))_{N \ge 1} \subset \mathcal{S}_{\mathbf{F}}^2([0,T]; \mathbb{R}^d)$:

THEOREM 4.1. *Let $\xi = (\xi^{(1)}, \xi^{(2)}) = \{((\xi_t^{(1,i)}(x))_{1 \le i \le d}, (\xi_t^{(2,i,j)}(x))_{1 \le i,j \le d}),$ $(t,x) \in [0,T] \times \mathbb{R}^d\}$ be a $d + d \times d$-dimensional continuous zero-mean Gaus-* sian field which is independent of the Brownian motion $W$ and has the co-variance functions

$$E[\xi_t^{(1,i)}(x)\xi_{t'}^{(1,j)}(x')] = \mathrm{cov}(b^i(x, X_t), b^j(x', X_{t'})),$$

$$E[\xi_t^{(1,i)}(x)\xi_{t'}^{(2,k,\ell)}(x')] = \mathrm{cov}(b^i(x, X_t), \sigma^{k,\ell}(x', X_{t'})),$$

$$E[\xi_t^{(2,i,j)}(x)\xi_{t'}^{(2,k,\ell)}(x')] = \mathrm{cov}(\sigma^{i,j}(x, X_t), \sigma^{k,\ell}(x', X_{t'}))$$

$$\text{for all } (t,x), (t', x') \in [0, T] \times \mathbb{R}^d, 1 \le i, j, k, \ell \le d,$$

and let $\overline{X} \in \mathcal{S}_{\mathbf{F}}^2([0, T]; \mathbb{R}^d)$ be the unique solution of the forward equation

$$
\begin{aligned}
(4.3) \quad \overline{X}_t = {}& \int_0^t \xi_s^{(1)}(X_s)\, ds + \int_0^t \xi_s^{(2)}(X_s)\, dW_s \\
& + \int_0^t (E[(\nabla_x b)(x, X_s)]_{x=X_s}\overline{X}_s + E[(\nabla_{x'} b)(x, X_s)\overline{X}_s]_{x=X_s})\, ds \\
& + \int_0^t (E[(\nabla_x \sigma)(x, X_s)]_{x=X_s}\overline{X}_s + E[(\nabla_{x'} \sigma)(x, X_s)\overline{X}_s]_{x=X_s})\, dW_s,
\end{aligned}
$$

$t \in [0, T]$, where $\overline{\mathbf{F}}$ is the filtration $\mathbf{F}$ augmented by $\sigma\{\xi_t(x), (t, x) \in [0, T] \times \mathbb{R}^d\}$. Then the sequence $(\sqrt{N}(X^N - X))_{N \ge 1}$ converges in law over $C([0, T]; \mathbb{R}^d)$ to the process $\overline{X}$.

REMARK 4.1. The continuity of the above introduced two-dimensional zero-mean Gaussian process $\xi$ is a direct consequence of Kolmogorov's Continuity Criterion for multi-parameter processes.

Indeed, a standard argument for mean-zero Gaussian random variables and a standard SDE estimate shows that, for all $m \ge 1$ and for some generic constant $C_m$ which can change from line to line,

$$
\begin{aligned}
E[|\xi_t(x) &- \xi_{t'}(x')|^{2m}] \\
&\le C_m \sum_{i=1}^d E[|\xi_t^{(1,i)}(x) - \xi_{t'}^{(1,i)}(x')|^{2m}] \\
&\quad + C_m \sum_{i,j=1}^d E[|\xi_t^{(2,i,j)}(x) - \xi_{t'}^{(2,i,j)}(x')|^{2m}] \\
&\le C_m \sum_{i=1}^d (E[|\xi_t^{(1,i)}(x) - \xi_{t'}^{(1,i)}(x')|^2])^m \\
&\quad + C_m \sum_{i,j=1}^d (E[|\xi_t^{(2,i,j)}(x) - \xi_{t'}^{(2,i,j)}(x')|^2])^m
\end{aligned}
$$

$$= C_m \sum_{i=1}^{d} (\text{Var}(b^i(x, X_t) - b^i(x', X_{t'})))^m$$

$$+ C_m \sum_{i,j=1}^{d} (\text{Var}(\sigma^{i,j}(x, X_t) - \sigma^{i,j}(x', X_{t'})))^m$$

$$\leq C_m(|t - t'|^m + |x - x'|^{2m}) \qquad \text{for all } (t, x), (t', x') \in [0, T] \times \mathbb{R}.$$

The proof of Theorem 4.1 will be split into a sequel of statements whose objective is the study of the limit behavior of the process $\sqrt{N}(X^N - X)$ as $N$ tends to infinity. The process $\sqrt{N}(X^N - X)$ can be described by the following SDE:

$$\sqrt{N}(X_t^N - X_t)$$

$$= \int_0^t \{\xi_s^{1,N}(X_s^N) + \sqrt{N}(E[b(x, X_s^N)]_{|x=X_s^N} - E[b(x, X_s)]_{|x=X_s})\} \, ds$$

$$+ \int_0^t \{\xi_s^{2,N}(X_s^N) + \sqrt{N}(E[\sigma(x, X_s^N)]_{|x=X_s^N} - E[\sigma(x, X_s)]_{|x=X_s})\} \, dW_s,$$

$$t \in [0, T],$$

where

$$\xi_t^{1,N}(x) = \frac{1}{\sqrt{N}} \sum_{k=1}^{N} (b(x, \Theta^k(X_t^N)) - E[b(x, X_t^N)]),$$

$$\xi_t^{2,N}(x) = \frac{1}{\sqrt{N}} \sum_{k=1}^{N} (\sigma(x, \Theta^k(X_t^N)) - E[\sigma(x, X_t^N)]),$$

$$(t, x) \in [0, T] \times \mathbb{R}^d.$$

Obviously, the random fields $\xi^{i,N} = \{\xi_t^{i,N}(x), (t, x) \in [0, T] \times \mathbb{R}^d\}$, $i = 1, 2$, are independent of the Brownian motion $W$ and they have their paths in $C([0, T] \times \mathbb{R}^d; \mathbb{R}^d)$ and $C([0, T] \times \mathbb{R}^d; \mathbb{R}^{d \times d})$, respectively. Moreover, we can characterize their limit behavior as follows:

PROPOSITION 4.1. *The sequence of the laws of the stochastic fields $\xi^N = (\xi^{1,N}, \xi^{2,N})$, $N \geq 1$, converges weakly on $C([0, T] \times \mathbb{R}^d; \mathbb{R}^d \times \mathbb{R}^{d \times d})$ toward the law of the continuous zero-mean Gaussian $(d + 1)$-parameter process $\xi = \{(\xi^{(1)}, \xi^{(2)})\}$ introduced in Theorem 4.1.*

The proof of the weak convergence of the laws $P \circ [\xi^N]^{-1}, N \geq 1$, on $C([0, T] \times \mathbb{R}^d; \mathbb{R}^d \times \mathbb{R}^{d \times d})$ is split into two lemmas: while the first lemma establishes the tightness of these laws on $C([0,T] \times \mathbb{R}^d; \mathbb{R}^d \times \mathbb{R}^{d \times d})$ the second lemma will study the convergence of their finite-dimensional marginal laws. The both results together, the tightness of the laws and the convergence of the finite-dimensional marginal laws imply the weak convergence.

LEMMA 4.1.   *The sequence of the laws of the stochastic processes* $\xi^N = (\xi^{1,N}, \xi^{2,N})$, $N \geq 1$, *is tight on* $C([0,T] \times \mathbb{R}^d; \mathbb{R}^d \times \mathbb{R}^{d \times d})$. *In fact, we have for all* $m \geq 1$ *the existence of some constant* $C_m$ *such that:*

(i) $E[|\xi_t^N(x)|^{2m}] \leq C_m$, *for all* $(t,x) \in [0,T] \times \mathbb{R}^d, N \geq 1$;

(ii) $E[|\xi_t^N(x) - \xi_{t'}^N(x')|^{2m}] \leq C_m(|t - t'|^m + |x - x'|^{2m})$, *for all* $(t,x), (t', x') \in [0,T] \times \mathbb{R}^d, N \geq 1$.

*Moreover, for all* $m \geq 1$ *there exists some constant* $C_m$ *such that*

(iii) $\quad E\left[\sup_{x \in \mathbb{R}^d} \left( \frac{|\xi_t^N(x)|}{1 + |x|} \right)^{2m}\right] \leq C_m \qquad$ *for all* $t \in [0,T], N \geq 1.$

PROOF.   We begin with the proof of statement (i). Putting

$$\overline{b}^{i,N}(t,x) := b^i(x, X_t^N) - E[b^i(x, X_t^N)]$$

and

$$\overline{\sigma}^{i,j,N}(t,x) := \sigma^{i,j}(x, X_t^N) - E[\sigma^{i,j}(x, X_t^N)], \qquad 1 \leq i, \ j \leq d, \ N \geq 1,$$

we have the existence of some generic constant $C_{d,m}$ such that, for all $(t,x) \in [0,T] \times \mathbb{R}^d$ and for all $N \geq 1$,

$$E[|\xi_t^N(x)|^{2m}] \leq C_{m,d}\left( \sum_{i=1}^d E\left[ \left| \frac{1}{\sqrt{N}} \sum_{k=1}^N \Theta^k(\overline{b}^{i,N}(t,x)) \right|^{2m} \right] \right.$$
$$\left. + \sum_{i,j=1}^d E\left[ \left| \frac{1}{\sqrt{N}} \sum_{k=1}^N \Theta^k(\overline{\sigma}^{i,j,N}(t,x)) \right|^{2m} \right] \right).$$

On the other hand, taking into account the independence of the zero-mean random variables $\Theta^k(\overline{b}^{i,N}(t,x)), 1 \leq k \leq N$ and the boundedness of the function $b$ we get, with the notation $\Gamma_{m,N} = \{(k_1, \dots, k_{2m}) \in \{1, \dots, N\}^{2m}$: for all $i(1 \leq i \leq 2m)$ there is a $j \in \{1, \dots, N\} \setminus \{i\}$ s.t. $k_i = k_j\}$,

$$E\left[ \left| \frac{1}{\sqrt{N}} \sum_{k=1}^N \Theta^k(\overline{b}^{i,N}(t,x)) \right|^{2m} \right]$$
$$= \frac{1}{N^m} \sum_{k_1, \dots, k_{2m}=1}^N E\left[ \prod_{j=1}^{2m} \Theta^{k_j}(\overline{b}^{i,N}(t,x)) \right]$$

$$= \frac{1}{N^m} \sum_{(k_1,\ldots,k_{2m}) \in \Gamma_{m,N}} E\left[\prod_{j=1}^{2m} \Theta^{k_j}(\overline{b}^{i,N}(t,x))\right]$$

$$\leq \frac{1}{N^m} \sum_{(k_1,\ldots,k_{2m}) \in \Gamma_{m,N}} \prod_{j=1}^{2m} E[|\Theta^{k_j}(\overline{b}^{i,N}(t,x))|^{2m}]^{1/(2m)}$$

$$\leq \frac{1}{N^m} \operatorname{card}(\Gamma_{m,N}) E[|\overline{b}^{i,N}(t,x)|^{2m}]$$

$$= C'_m E[|\overline{b}^{i,N}(t,x)|^{2m}] \leq C_m$$

for some $C'_m,\ C_m \in \mathbb{R}_+$, which are independent of $N \geq 1,\ 1 \leq i \leq d$ and $(t,x) \in [0,T] \times \mathbb{R}^d$. Similarly, we see that, for all $N \geq 1,\ 1 \leq i,j \leq d$ and $(t,x) \in [0,T] \times \mathbb{R}^d$,

$$E\left[\left|\frac{1}{\sqrt{N}} \sum_{k=1}^N \Theta^k(\overline{\sigma}^{i,j,N}(t,x))\right|^{2m}\right] \leq C_m.$$

Consequently, for some $C_m \in \mathbb{R}_+$ neither depending on $N \geq 1$ nor on $(t,x) \in [0,T] \times \mathbb{R}^d$, $E[|\xi_t^N(x)|^{2m}] \leq C_m$, and the same argument, but now applied to $\nabla_x b(x,x'), \nabla_x \sigma(x,x')$, also yields $E[|\nabla_x \xi_t^N(x)|^{2m}] \leq C_m$.

For proving statement (iii) of the lemma we observe that, for all $m \geq (d+1)/2$, due to Morrey's inequality and the estimates obtained above, we have the existence of some generic constant $C_m$ depending on $m$ and $d$ such that, with the notation $\widehat{\xi}_t^N(x) = (1+|x|)^{-1}\xi_t^N(x)$ and $\gamma = 1 - \frac{d}{2m}$,

$$E\left[\sup_{x \in \mathbb{R}^d}\left(\frac{|\xi_t^N(x)|}{1+|x|}\right)^{2m}\right]$$

$$\leq E\left[\left(\sup_{x \in \mathbb{R}^d}|\widehat{\xi}_t^N(x)| + \sup_{x \neq x'}\frac{|\widehat{\xi}_t(x) - \widehat{\xi}_t(x')|}{|x-x'|^\gamma}\right)^{2m}\right]$$

$$\leq C_m E\left[\int_{\mathbb{R}^d}(|\widehat{\xi}_t^N(x)|^{2m} + |\nabla_x \widehat{\xi}_t^N(x)|^{2m})\,dx\right]$$

$$\leq C_m \int_{\mathbb{R}^d}(E[|\xi_t^N(x)|^{2m}] + E[|\nabla_x \xi_t^N(x)|^{2m}])\frac{dx}{(1+|x|)^{2m}}$$

$$\leq C_m$$

for all $t \in [0,T], N \geq 1$. To get (ii) we follow the argument given for the proof of (i). Combining it with a standard SDE estimate we get, for a generic constant $C_m$ and all $(t,x),(t',x') \in [0,T] \times \mathbb{R}^d$ and $N \geq 1$,

$$E[|\xi_t^N(x) - \xi_{t'}^N(x')|^{2m}]$$

$$\leq C_m\left(\sum_{i=1}^d E[|\overline{b}^{i,N}(t,x) - \overline{b}^{i,N}(t',x')|^{2m}]\right.$$

$$+ \sum_{i,j=1}^{d} E[|\overline{\sigma}^{i,j,N}(t,x) - \overline{\sigma}^{i,j,N}(t',x')|^{2m}] \Bigg)$$

$$\leq C_m(E[|b(x,X_t^N) - b(x',X_{t'}^N)|^{2m}] + E[|\sigma(x,X_t^N) - \sigma(x',X_{t'}^N)|^{2m}])$$

$$\leq C_m(|t - t'|^m + |x - x'|^{2m}).$$

Finally, to complete the proof it only remains to observe that due to Kolmogorov's weak compactness criterion for multi-parameter processes the estimates (i) and (ii) imply the tightness of the sequence of laws of $\xi^N$, $N \geq 1$, on $C([0,T] \times \mathbb{R}^d; \mathbb{R}^d \times \mathbb{R}^{d \times d})$. □

For proving Proposition 4.1 we have still to show the following lemma.

LEMMA 4.2. *The finite-dimensional laws of the stochastic fields $\xi^N = (\xi^{1,N}, \xi^{2,N})$, $N \geq 1$, converge weakly to the corresponding finite-dimensional laws of the continuous zero-mean Gaussian $(d+1)$-parameter process $\xi = \{(\xi_t^{(1)}(x), \xi_t^{(2)}(x)), (t,x) \in [0,T] \times \mathbb{R}^d\}$ introduced in Theorem 4.1.*

PROOF. We decompose the random field $\xi^N = (\xi^{1,N}, \xi^{2,N})$. For this we first consider the component $\xi^{2,N}$ and we represent it as the sum of the random fields $\xi^{2,N,1}$ and $\xi^{2,N,2}$, where

$$\xi_t^{2,N,1}(x) = \frac{1}{\sqrt{N}} \sum_{k=1}^{N} (\Theta^k(\sigma(x,X_t)) - E[\sigma(x,X_t)]), \qquad (t,x) \in [0,T] \times \mathbb{R}^d,$$

$$\xi_t^{2,N,2}(x) = \frac{1}{\sqrt{N}} \sum_{k=1}^{N} (\Theta^k(\sigma(x,X_t^N) - \sigma(x,X_t)) - E[\sigma(x,X_t^N) - \sigma(x,X_t)])$$

for $(t,x) \in [0,T] \times \mathbb{R}^d$. Since the bounded $(d+1)$-parameter fields $\Theta^k(\sigma(\cdot,X)) = \{\sigma(x,\Theta^k(X_t)), (t,x) \in [0,T] \times \mathbb{R}^d\}$, $k \geq 0$, are i.i.d., it follows directly from the central limit theorem that the finite-dimensional laws of the random field $\xi^{2,N,1}$ converge weakly to the corresponding finite-dimensional laws of the zero-mean Gaussian $(d+1)$-parameter process $\xi^{(2)}$ whose covariance function coincides with that of the field $\sigma(\cdot,X) = \{\sigma(x,X_t), (t,x) \in [0,T] \times \mathbb{R}^d\}$ [recall the definition of $\xi = (\xi^{(1)}, \xi^{(2)})$ given in Theorem 4.1]. On the other hand, since due to Proposition 3.2,

$$E[|\xi_t^{2,N,2}(x)|^2]$$

$$= \frac{1}{N} \sum_{i,j=1}^{d} E\Bigg[\Bigg|\sum_{k=1}^{N} (\Theta^k(\sigma^{i,j}(x,X_t^N) - \sigma^{i,j}(x,X_t))$$

$$- E[\sigma^{i,j}(x,X_t^N) - \sigma^{i,j}(x,X_t)])\Bigg|^2\Bigg]$$

$$= \sum_{i,j=1}^{d} \mathrm{Var}(\sigma^{i,j}(x, X_t^N) - \sigma^{i,j}(x, X_t)) \leq CE[|X_t^N - X_t|^2]$$

$$\leq \frac{C}{N} \qquad \text{for all } N \geq 1, (t,x) \in [0,T] \times \mathbb{R}^d,$$

it follows that also the finite-dimensional laws of the random field $\xi^{2,N} = \xi^{2,N,1} + \xi^{2,N,2}$ converge weakly to the corresponding finite-dimensional laws of the zero-mean Gaussian field $\xi^{(2)}$.

By applying now the argument developed above to the couple $\xi^N = (\xi^{1,N}, \xi^{2,N})$ we can complete the proof.  □

As already mentioned before, Proposition 4.1 follows directly from the Lemmas 4.1 and 4.2. Proposition 4.1 again allows to apply Skorohod's representation theorem. Taking into account that the random fields $\xi^N$, $N \geq 1$, are all independent of the Brownian motion $W$, we can conclude from Skorohod's representation theorem that, on an appropriate complete probability space $(\Omega', \mathcal{F}', P')$ there exist $(d + d \times d)$-dimensional $(d+1)$-parameter processes $\xi'^{tN} = \{\xi_t'^N(x), (t,x) \in [0,T] \times \mathbb{R}^d\}$ $(N \geq 1)$ and $\xi' = \{\xi_t'(x), (t,x) \in [0,T] \times \mathbb{R}^d\}$, as well as a $d$-dimensional Brownian motion $W' = (W_t')_{t \in [0,T]}$, such that:

(i)   $P_{\xi'}' = P_\xi, P_{\xi'^{tN}}' = P_{\xi^N}, N \geq 1$;
(ii)  $\xi'^{tN} = (\xi'^{1,N}, \xi'^{2,N}) \longrightarrow \xi' = (\xi'^1, \xi'^2)$, uniformly on compacts, $P'$-a.s.;
(iii) $W'$ is independent of $\xi'$ and $\xi'^{tN}$, for all $N \geq 1$.

For this new probability space we introduce the filtration $\mathbf{F}' = (\mathcal{F}_t' = \mathcal{F}_t^{W'} \vee \mathcal{F}_0')_{t \in [0,T]}$, where $\mathcal{F}_0' = \sigma\{\xi_s'(x), \zeta_s'^N(x), (s,x) \in [0,T] \times \mathbb{R}^d, N \geq 1\} \vee \mathcal{N}$, and we observe that $W'$ is an $\mathbf{F}'$-Brownian motion.

Given $\xi'^{tN}, \xi'$ and the Brownian motion $W'$ we now redefine the processes $X^N$ and $X$ on our new probability space. For this we recall that the process $X^N$ is defined as unique solution of (4.1). In virtue of the definition of $\xi^N = (\xi^{1,N}, \xi^{2,N})$ we can rewrite (4.1) as follows:

$$X_t^N = x_0 + \int_0^t \left( \frac{1}{\sqrt{N}} \xi_s^{1,N}(X_s^N) + E[b(x, X_s^N)]_{|x=X_s^N} \right) ds$$

$$+ \int_0^t \left( \frac{1}{\sqrt{N}} \xi_s^{2,N}(X_s^N) + E[\sigma(x, X_s^N)]_{|x=X_s^N} \right) dW_s, \qquad t \in [0,T].$$

Taking into account that, $P$-a.s.,

$$|\xi_s^{i,N}(x) - \xi_s^{i,N}(x')| \leq 2\sqrt{N} L|x - x'| \qquad \text{for all } s \in [0,T], x, x' \in \mathbb{R}^d, i = 1, 2,$$

where $L$ denotes the Lipschitz constant of $b$ and of $\sigma$, the fact that $\xi^N = (\xi^{1,N}, \xi^{2,N})$ and $\xi'^{tN} = (\xi'^{1,N}, \xi'^{2,N})$ obey the same law has as consequence that also $\xi_s'^{1,N}(\cdot)$ and $\xi_s'^{2,N}(\cdot)$ are Lipschitz, with the same Lipschitz constant as $\xi_s^{1,N}$ and $\xi_s^{2,N}$. Therefore, the equation

$$X'^N_t = x_0 + \int_0^t \left( \frac{1}{\sqrt{N}} \xi_s'^{1,N}(X'^N_s) + E'[b(x, X'^N_s)]_{|x=X'^N_s} \right) ds$$

$$+ \int_0^t \left( \frac{1}{\sqrt{N}} \xi_s'^{2,N}(X'^N_s) + E'[\sigma(x, X'^N_s)]_{|x=X'^N_s} \right) dW_s, \qquad t \in [0, T],$$

admits a unique solution process $X'^N \in \mathcal{S}^2_{\mathbf{F}'}([0, T]; \mathbb{R}^d)$. In the same spirit in which we have associated with $X^N$ the process $X'^N$ we define the analogue to $X$ on the probability space $(\Omega', \mathcal{F}', P')$ by letting $X' \in \mathcal{S}^2_{\mathbf{F}'}([0, T]; \mathbb{R}^d)$ be the unique solution of the SDE

$$(4.4) \quad X'_t = x_0 + \int_0^t E'[b(x, X'_s)]_{|x=X'_s} ds + \int_0^t E'[\sigma(x, X'_s)]_{|x=X'_s} dW'_s,$$
$$t \in [0, T].$$

We remark that the solution $X'$ is adapted to the filtration generated by the Brownian motion $W'$. Moreover, the fact that the laws of $(\xi'^N, W')$ and $(\xi^N, W)$ coincide implies that also

$$P' \circ (X', X'^N, \xi'^N, W')^{-1} = P \circ (X, X^N, \xi^N, W)^{-1}, \qquad N \geq 1.$$

Consequently, Theorem 4.1 follows immediately from the following proposition:

PROPOSITION 4.2. Let $\widehat{X} \in \mathcal{S}^2_{\mathbf{F}'}([0, T]; \mathbb{R}^d)$ be the unique solution of the SDE

$$(4.5) \quad \begin{aligned} \widehat{X}_t &= \int_0^t (\xi_s'^1(X'_s) + E'[\nabla_x b(x, X'_s)]_{|x=X'_s} \widehat{X}_s \\ &\qquad + E'[\nabla_{x'} b(x, X'_s) \widehat{X}_s]_{|x=X'_s}) ds \\ &\qquad + \int_0^t (\xi_s'^2(X'_s) + E'[\nabla_x \sigma(x, X'_s)]_{|x=X'_s} \widehat{X}_s \\ &\qquad + E'[\nabla_{x'} \sigma(x, X'_s) \widehat{X}_s]_{|x=X'_s}) dW'_s, \end{aligned}$$

$t \in [0, T]$. Then

$$E'\left[ \sup_{s \in [0, t]} |\sqrt{N}(X'^N_s - X'_s) - \widehat{X}_s|^2 \right] \longrightarrow 0 \qquad as \ N \to +\infty.$$

PROOF. We begin by showing

STEP 1. The existence and uniqueness of the solution $\widehat{X} \in \mathcal{S}_{\mathbf{F}'}([0, T]; \mathbb{R}^d)$ of the above SDE (4.5).

For this we remark that, thanks to the independence of $\xi' = (\xi'^1, \xi'^2)$ and $W'$ we have also that of $\xi' = (\xi'^1, \xi'^2)$ and $X'$, so that

$$E'[|\xi_s'^i(X_s')|^{2m}] = \int_{\mathbb{R}^d} E'[|\xi_s'^i(x)|^{2m}] P'_{X_s'}(dx) \le C_m, \qquad s \in [0, T], \ m \ge 1$$

(cf. Lemma 4.1). Since the other coefficients of the above linear SDE are bounded we have the existence and the uniqueness of the solution $\widehat{X}$ in $\mathcal{S}_{\mathbf{F}'}([0, T]; \mathbb{R}^d)$; moreover, $\widehat{X}$ is adapted with respect to the filtration generated by $W'$ and $E'[\sup_{t \in [0,T]} |\widehat{X}_t|^p] < \infty$, for all $p \ge 1$.

STEP 2. $(\xi_t'^N(X_t'^N))_{t \in [0,T]} \longrightarrow (\xi_t'(X_t'))_{t \in [0,T]}$ in $L^2([0, T] \times \Omega', dt \, dP')$, as $N \to +\infty$.

Indeed, for all $m \ge 1$, we have

$$E'[|\xi_s'^N(X_s'^N)|^m] = E[|\xi_s^N(X_s^N)|^m]$$
$$\le \left( E\left[ \sup_{x \in \mathbb{R}^d} \left( \frac{|\xi_s^N(x)|}{1 + |x|} \right)^{2m} \right] \right)^{1/2} (E[(1 + |X_s^N|)^{2m}])^{1/2}.$$

Thus, using Lemma 4.1(iii) and the fact that the coefficients of (4.1) are bounded, uniformly with respect to $N \ge 1$, we can conclude that, for some constant $C_m$,

$$E'[|\xi_s'^N(X_s'^N)|^m] \le C_m \qquad \text{for all } s \in [0, T], N \ge 1.$$

On the other hand, recalling that $\xi'^N = (\xi'^{1,N}, \xi'^{2,N}) \to \xi' = (\xi'^1, \xi'^2)$ uniformly on compacts of $[0, T] \times \mathbb{R}^d$, $P'$-a.s., and

$$E'\left[ \sup_{s \in [0,T]} |X_s'^N - X_s'|^2 \right] = E\left[ \sup_{s \in [0,T]} |X_s^N - X_s|^2 \right] \le \frac{C}{N},$$

$$N \ge 1 \text{ (cf. Proposition 3.2)},$$

we see that $\sup_{s \in [0,T]} |\xi_s'^N(X_s'^N) - \xi_s'(X_s')| \to 0$, in probability. Combining this with the uniform square integrability of $\{\xi_s'^N(X_s'^N), s \in [0, T], N \ge 1\}$, proved above, yields the convergence of $\xi'^N(X'^N)$ to $\xi'(X')$ in $L^2([0, T] \times \Omega', dt \, dP')$.

STEP 3. Let us now decompose $\sqrt{N}(X_t'^N - X_t') - \widehat{X}_t$ as follows:

$$\sqrt{N}(X_t'^N - X_t') - \widehat{X}_t = I_t^{1,N} + I_t^{2,N} + I_t^{3,N}, \qquad t \in [0, T],$$

where

$$I_t^{1,N} = \int_0^t (\xi_s'^{1,N}(X_s'^N) - \xi_s'^1(X_s')) \, ds + \int_0^t (\xi_s'^{2,N}(X_s'^N) - \xi_s'^2(X_s')) \, dW_s',$$

$$I_t^{2,N} = \int_0^t \{\sqrt{N}(E'[b(x, X_s'^N)]|_{x=X_s'^N} - E'[b(x, X_s')]|_{x=X_s'})$$

$$- (E'[\nabla_x b(x, X'_s)]_{|x=X'_s} \widehat{X}_s + E'[\nabla_{x'} b(x, X'_s) \widehat{X}_s]_{|x=X'_s}) \} \, ds,$$

$$I_t^{3,N} = \int_0^t \{ \sqrt{N} (E'[\sigma(x, X'^N_s)]_{|x=X'^N_s} - E'[\sigma(x, X'_s)]_{|x=X'_s})$$

$$- (E'[\nabla_x \sigma(x, X'_s)]_{|x=X'_s} \widehat{X}_s + E'[\nabla_{x'} \sigma(x, X'_s) \widehat{X}_s]_{|x=X'_s}) \} \, dW'_s.$$

From the Step 2 we know already that $E'[\sup_{t \in [0,T]} |I_t^{1,N}|^2] \to 0$, as $N \to +\infty$. Let us now estimate $E'[\sup_{t \in [0,T]} |I_t^{3,N}|^2]$. For this we remark that, thanks to the continuous differentiability of the function $\sigma$, with the notation $X'^{N,\gamma}_s := X'_s + \gamma(X'^N_s - X'_s)$ we get

$$I_t^{3,N} = \int_0^t \left\{ \int_0^1 E'[\nabla_x \sigma(x_1, X'^{N,\gamma}_s) - \nabla_x \sigma(x_2, X'_s)]_{|x_1=X'^{N,\gamma}_s, x_2=X'_s} \, d\gamma \times \widehat{X}_s \right.$$

$$+ \int_0^1 E'[\{ \nabla_{x'} \sigma(x_1, X'^{N,\gamma}_s)$$

$$- \nabla_{x'} \sigma(x_2, X'_s) \} \widehat{X}_s]_{|x_1=X'^{N,\gamma}_s, x_2=X'_s} \, d\gamma \bigg\} \, dW'_s$$

$$+ \int_0^t \{ E'[\nabla_x \sigma(x, X'^{N,\gamma}_s)]_{x=X'^{N,\gamma}_s} (\sqrt{N}(X'^N_s - X'_s) - \widehat{X}_s)$$

$$+ E'[\nabla_{x'} \sigma(x, X'^{N,\gamma}_s) (\sqrt{N}(X'^N_s - X'_s) - \widehat{X}_s)]_{x=X'^{N,\gamma}_s} \} \, dW'_s,$$
$$t \in [0,T],$$

and taking into account that the first-order derivatives of $\sigma$ are bounded and continuous, it follows from the convergence of $X'^{N,\gamma} \to X'$ ($N \to +\infty$) in $\mathcal{S}^2_{\mathbf{F}'}([0,T]; \mathbb{R}^d)$, for all $\gamma \in [0,1]$, that for some sequence $0 < \rho_N \searrow 0$,

$$E'\left[ \sup_{s \in [0,t]} |I_s^{3,N}|^2 \right] \le \rho_N + C \int_0^t E'[|\sqrt{N}(X'^N_s - X'_s) - \widehat{X}_s|^2] \, ds,$$
$$t \in [0,T], \ N \ge 1.$$

The same argument yields the above estimate also for $I^{2,N}$. Consequently,

$$E'\left[ \sup_{s \in [0,t]} |\sqrt{N}(X'^N_s - X'_s) - \widehat{X}_s|^2 \right]$$

$$\le E'\left[ \sup_{t \in [0,T]} |I_t^{1,N}|^2 \right] + 2\rho_N + 2C \int_0^t E'[|\sqrt{N}(X'^N_s - X'_s) - \widehat{X}_s|^2] \, ds,$$
$$t \in [0,T], \ N \ge 1,$$

and from Gronwall's lemma we obtain

$$E'\left[ \sup_{s \in [0,t]} |\sqrt{N}(X'^N_s - X'_s) - \widehat{X}_s|^2 \right] \longrightarrow 0 \qquad \text{as } N \to +\infty.$$

The proof is complete. □

It remains to give the proof of Theorem 4.1.

PROOF OF THEOREM 4.1. In analogy to the existence and uniqueness of the solution $\widehat{X} \in \mathcal{S}^2_{\mathbf{F}'}([0,T];\mathbb{R}^d)$ we can prove that of $\overline{X} \in \mathcal{S}^2_{\mathbf{F}}([0,T];\mathbb{R}^d)$. Moreover, since $W$ is independent of $\xi$ and $\xi^N, N \geq 1$, and, on the other hand, also $W'$ is independent of $\xi'$ and $\xi'^N, N \geq 1$, it follows that:

(i) $P \circ (W, \xi^N)^{-1} = P' \circ (W', \xi'^N)^{-1}$, for all $N \geq 1$;

(ii) $P \circ (W, \xi)^{-1} = P' \circ (W', \xi')^{-1}$. But $X'^N$ is the unique solution of an SDE governed by the $W'$ and $\xi'^N$, and $X^N$ is the unique solution of the same equation but driven by $W$ and $\xi^N$ instead of $W'$ and $\xi'^N$. On the other hand, $X'$ is the unique solution of (4.2) with $W'$ at the place of $W$. Thus,

(iii) $P \circ (W, \xi^N, X, X^N)^{-1} = P' \circ (W', \xi'^N, X', X'^N)^{-1}$, for all $N \geq 1$. Concerning the process $\widehat{X}$, it is the unique solution of the same equation as $\overline{X}$, but with $\xi$ and $W$ replaced by $\xi'$ and $W'$, respectively. This, together with (ii), yields

(iv) $P \circ (W, \xi, X, \overline{X})^{-1} = P' \circ (W', \xi', X', \widehat{X})^{-1}$. From (iii) and (iv) we see, in particular, that $P_{\sqrt{N}(X^N - X)} = P'_{\sqrt{N}(X'^N - X')}$, $N \geq 1$, and $P_{\overline{X}} = P'_{\widehat{X}}$. For completing the proof it suffices now to remark that the convergence of $\sqrt{N}(X'^N - X')$ to $\widehat{X}$ in $\mathcal{S}^2_{\mathbf{F}'}([0,T];\mathbb{R}^d)$ implies that of their laws on $C([0,T];\mathbb{R}^d)$. □

4.2. *Limit behavior of $\sqrt{N}(X^N - X, Y^N - Y, Z^N - Z)$.* As in the preceding section we suppose that $b:\mathbb{R}^d \times \mathbb{R}^d \to \mathbb{R}^d$ and $\sigma:\mathbb{R}^d \times \mathbb{R}^d \to \mathbb{R}^{d \times d}$ are bounded, continuously differentiable functions with bounded first-order derivatives, and $x_0$ is an arbitrarily fixed element of $\mathbb{R}^d$. With the forward equation

(4.6)
$$X^N_t = x_0 + \int_0^t \frac{1}{N} \sum_{k=1}^N b(X^N_s, \Theta^k(X^N_s)) \, ds$$
$$+ \int_0^t \frac{1}{N} \sum_{k=1}^N \sigma(X^N_s, \Theta^k(X^N_s)) \, dW_s, \qquad t \in [0,T],$$

we associate the BSDE

(4.7)
$$Y^N_t = \frac{1}{N} \sum_{k=1}^N \Phi(X^N_T, \Theta^k(X^N_T)) + \int_t^T \frac{1}{N} \sum_{k=1}^N f(\Lambda^N_s, \Theta^k(\Lambda^N_s)) \, ds$$
$$- \int_t^T Z^N_s \, dW_s,$$

where

$$\Lambda^N = (X^N, Y^N, Z^N) \in B_{\mathbf{F}}^2 := \mathcal{S}_{\mathbf{F}}^2([0,T]; \mathbb{R}^d) \times \mathcal{S}_{\mathbf{F}}^2([0,T]; \mathbb{R}) \times L_{\mathbf{F}}^2([0,T]; \mathbb{R}^d)$$

is composed of the unique solution $X^N$ of the forward SDE and the unique solution $(Y^N, Z^N)$ of the BSDE. The functions $f : (\mathbb{R}^d \times \mathbb{R} \times \mathbb{R}^d)^2 \to \mathbb{R}$ and $\Phi : \mathbb{R}^d \times \mathbb{R}^d \to \mathbb{R}$ are supposed to be bounded and continuously differentiable with bounded first-order derivatives. Moreover, our approach imposes the following additional hypothesis for the function $f$:

(E) $f((x,y,z),(x',y',z')) = f((x,y,z),(x',y'))$, *for all* $(x,y,z)$, $(x,'y',z')$ $\in \mathbb{R}^d \times \mathbb{R} \times \mathbb{R}^d$.

We know already that

$$E\left[\sup_{t \in [0,T]} (|X_t^N - X_t|^2 + |Y_t^N - Y_t|^2) + \int_0^T |Z_t^N - Z_t|^2\right] \le \frac{C}{N}, \qquad N \ge 1,$$

for some constant $C$, and we are interested in the description of the limit behavior of the sequence $\sqrt{N}(X^N - X, Y^N - Y, Z^N - Z)$. A crucial role in these studies will be played by the following uniform estimate for $Z_t^N, t \in [0,T], N \ge 1$.

LEMMA 4.3. *Under the above assumptions on the coefficients there exists some real constant $C$ such that, $dt\,dP$-a.e.,*

$$|Z_t| \le C \quad and \quad |Z_t^N| \le C \quad for\ all\ N \ge 1.$$

PROOF. We remark that $\Lambda^N = (X^N, Y^N, Z^N) \in B_{\mathbf{F}}^2$ is the unique solution of the system

$$X_t^N = x_0 + \int_0^t b_s^N(X_s^N)\,ds + \int_0^t \sigma_s^N(X_s^N)\,dW_s,$$

$$Y_t^N = \Phi^N(X_T^N) + \int_t^T f_s^N(\Lambda_s^N)\,ds - \int_t^T Z_s^N\,dW_s, \qquad t \in [0,T],$$

where the coefficients

$$\sigma_t^N(x) = \frac{1}{N}\sum_{k=1}^N \sigma(x, \Theta^k(X_t^N)); \qquad b_t^N(x) = \frac{1}{N}\sum_{k=1}^N b(x, \Theta^k(X_t^N)),$$

$$\Phi^N(x) = \frac{1}{N}\sum_{k=1}^N \Phi(x, \Theta^k(X_T^N)); \qquad f_t^N(x,y,z) = \frac{1}{N}\sum_{k=1}^N f((x,y,z), \Theta^k(\Lambda_t^N)),$$

$t \in [0,T], (x,y,z) \in \mathbb{R}^d \times \mathbb{R} \times \mathbb{R}^d$, are bounded random fields such that:

(i) $\sigma_t^N(x), b_t^N(x), \Phi^N(x)$ and $f_t^N(x,y,z)$ are $\mathcal{F}_0$-measurable, for all $(t,x,y,z)$, and, hence, independent of the driving Brownian motion $W$;

(ii) $\sigma_t^N(\cdot), b_t^N(\cdot), \Phi^N(\cdot)$ and $f_t^N(\cdot)$ are continuously differentiable and their derivatives are bounded by some constant which does neither depend on $t \in [0,T]$ nor on $\omega \in \Omega$.

Let $D_t = (D_t^1, \ldots, D_t^d), t \in [0,T]$, denote the Malliavin derivative with respect to the Brownian motion $W = (W^1, \ldots, W^d)$. For a smooth functional $F \in \mathcal{S}$ of the form $F = \zeta \varphi(W_{t_1}^{i_1}, \ldots, W_{t_n}^{i_n})$, with $n \geq 1, t_1, \ldots, t_n \in [0,T], 1 \leq i_1, \ldots, i_n \leq n, \varphi \in C_b^\infty(R^{n \times d})$ and $\zeta \in L^\infty(\Omega, \mathcal{F}_0, P)$, the Malliavin derivative $D_t^j F$ is defined by

$$D_t^j F = \zeta \sum_{\ell=1}^n \partial_{x_\ell} \varphi(W_{t_1}^{i_1}, \ldots, W_{t_n}^{i_n}) I_{[0,t_\ell]}(t) \delta_{i_\ell, j}, \qquad t \in [0,T], 1 \leq j \leq d.$$

The operator $D : \mathcal{S} \subset L^2(\Omega, \mathcal{F}, P) \to L^2(\Omega \times [0,T], dP\,dt; R^d)$ is closable; its closure is denoted by $(D, \mathbb{D}^{1,2})$.

It is well known that, for the coefficients $b^N, \sigma^N, f^N$ and $\Phi^N$ with the above properties, $X^N, Y^N$ and $Z^N$ belong to $L^2([0,T]; \mathbb{D}^{1,2}), X_t^N, Y_t^N \in \mathbb{D}^{1,2}$, for all $t \in [0,T]$, and

$$D_s^i X_t^N = \sigma_s^{N,i,\cdot}(X_s^N) + \int_s^t D_s^i X_r^N \nabla_x b_r^N(X_r^N)\,dr + \int_s^t D_s^i X_r^N \nabla_x \sigma_r^N(X_r^N)\,dW_r,$$

$$D_s^i Y_t^N = D_s^i X_T^N \nabla_x \Phi^N(X_T^N) + \int_t^T D_s^i \Lambda_r^N \nabla_\lambda f_r^N(\Lambda_r^N)\,dr - \int_t^T D_s^i Z_r^N\,dW_r,$$

$0 \leq s \leq t \leq T, 1 \leq i \leq d$, where $\sigma_s^{N,i,\cdot}$ denotes the $i$th row of the matrix $\sigma_s^N$. The above result is standard in the Malliavin calculus, the interested reader is referred, for example, to the book of Nualart [9]. From standard SDE estimates we then get

$$E\left[\sup_{t \in [s,T]} |D_s Y_t^N|^2 \big| \mathcal{F}_s\right] \leq C, \qquad P\text{-a.s., for all } t \in [0,T], N \geq 1.$$

On the other hand, differentiating the equation

$$Y_s^N = Y_0^N - \int_0^s f_r^N(\Lambda_r^N)\,dr + \int_0^s Z_r^N\,dW_r, \qquad s \in [0,T],$$

in Malliavin's sense we obtain, for $0 \leq s \leq t \leq T$ and $1 \leq i \leq d$,

$$D_s^i Y_t^N = Z_s^{N,i} - \int_s^t D_s^i \Lambda_r^N f_r^N(\Lambda_r^N)\,dr + \int_s^t D_s^i Z_r^N\,dW_r.$$

Consequently, letting $s \searrow t$ we have $D_s Y_t^N \to Z_s^N$ in $L^2(\Omega, \mathcal{F}, P)$, $ds$-a.e. From there and the above estimate for $D_s^i Y_t^N$ we deduce that $|Z_s^N| \leq C$, $ds\,dP$-a.e., for all $N \geq 1$. Finally, since this estimate is uniform with respect to $N$, it holds also true for the limit $Z$ of the sequence $Z^N, N \geq 1$. The proof is complete. $\square$

From the above lemma and the assumption of boundedness of the coefficients $\Phi$ and $f$ of the BSDEs we get the following by a standard estimate:

COROLLARY 4.1. *For all $m \geq 1$ there is some constant $C_m$ such that, for all $t$, $t' \in [0,T]$ and all $N \geq 1$,*

$$E[|Y_t^N - Y_{t'}^N|^{2m}] \leq C_m |t - t'|^m.$$

*The same estimate holds true for the limit $Y$ of the sequence $Y^N, N \geq 1$.*

Recall that the limit $X$ of the sequence $X^N, N \geq 1$, is the unique solution of the forward equation

$$
\begin{aligned}
(4.8) \quad X_t = x_0 + \int_0^t E[b(x, X_s)]_{|x=X_s}\, ds + \int_0^t E[\sigma(x, X_s)]_{|x=X_s}\, dW_s, \\
t \in [0, T],
\end{aligned}
$$

while the limit of the sequence $(Y^N, Z^N)$ is given by the unique solution $(Y, Z)$ of the BSDE

$$
\begin{aligned}
(4.9) \quad Y_t = E[\Phi(x, X_T)]_{x=X_T} + \int_t^T E[f(\lambda, \Lambda_s)]_{\lambda=\Lambda_s}\, ds - \int_t^T Z_s\, dW_s, \\
t \in [0, T],
\end{aligned}
$$

where $\Lambda = (X, Y, Z) \in B_{\mathbf{F}}^2$.

We have the following limit behavior of $(\sqrt{N}(X^N - X, Y^N - Y, Z^N - Z))_{N \geq 1} \subset B_{\mathbf{F}}^2$.

THEOREM 4.2. *Let $\xi = (\xi^{(1)}, \xi^{(2)}, \xi^{(3)}, \xi^{(4)}) = \{((\xi_t^{(1,i)}(x))_{1 \leq i \leq d},$ $(\xi_t^{(2,i,j)}(x))_{1 \leq i,j \leq d}, \xi^{(3)}(x), \xi_t^{(4)}(x, y, z)), (t, x, y, z) \in [0, T] \times \mathbb{R}^d \times \mathbb{R} \times \mathbb{R}^d\}$ be a $(d + d \times d + 2)$-dimensional continuous zero-mean Gaussian field which is independent of the Brownian motion $W$ and has the covariance function*

$$
E\left[
\begin{pmatrix}
\xi_t^{(1)}(x) \\
\xi_t^{(2)}(x) \\
\xi^{(3)}(x) \\
\xi_t^{(4)}(x, y, z)
\end{pmatrix}
\otimes
\begin{pmatrix}
\xi_{t'}^{(1)}(x') \\
\xi_{t'}^{(2)}(x') \\
\xi^{(3)}(x') \\
\xi_{t'}^{(4)}(x', y', z')
\end{pmatrix}
\right]
$$

$$
= E\left[
\begin{pmatrix}
b(x, X_t) - E[b(x, X_t)] \\
\sigma(x, X_t) - E[\sigma(x, X_t)] \\
\Phi(x, X_T) - E[\Phi(x, X_T)] \\
[f(\lambda, t)]
\end{pmatrix}
\otimes
\begin{pmatrix}
b(x', X_{t'}) - E[b(x', X_{t'})] \\
\sigma(x', X_{t'}) - E[\sigma(x', X_{t'})] \\
\Phi(x', X_T) - E[\Phi(x', X_T)] \\
[f(\lambda', t')]
\end{pmatrix}
\right],
$$

*where $[f(\lambda, t)] := f(\lambda, \Lambda_t) - E[f(\lambda, \Lambda_t)]$, $\lambda = (x, y, z), \lambda' = (x', y', z')$, $(t, x, y, z), (t', x', y', z') \in [0, T] \times \mathbb{R}^d \times \mathbb{R} \times \mathbb{R}^d$.*

*Recall that, for $a = (a^i)_{1 \le i \le m}, b = (b^i)_{1 \le i \le m} \in \mathbb{R}^m$, $a \otimes b = (a^i b^j)_{1 \le i, j \le m} \in \mathbb{R}^{m \times m}$; $\xi_t^{(2)}(x) = (\xi_t^{(2,i,j)}(x))_{1 \le i, j \le d}$ and $\sigma(x, x') = (\sigma^{i,j}(x, x'))_{1 \le i, j \le d}$ are regarded as $d^2$-dimensional vectors. By $\overline{\Lambda} = (\overline{X}, \overline{Y}, \overline{Z}) \in B_{\overline{\mathbf{F}}}^2$ we denote the unique solution of the system*

$$\overline{X}_t = \int_0^t \xi_s^{(1)}(X_s) \, ds + \int_0^t \xi_s^{(2)}(X_s) \, dW_s$$

$$\text{(4.10)} \qquad + \int_0^t (E[(\nabla_x b)(x, X_s)]_{x = X_s} \overline{X}_s + E[(\nabla_{x'} b)(x, X_s) \overline{X}_s]_{x = X_s}) \, ds$$

$$+ \int_0^t (E[(\nabla_x \sigma)(x, X_s)]_{x = X_s} \overline{X}_s + E[(\nabla_{x'} \sigma)(x, X_s) \overline{X}_s]_{x = X_s}) \, dW_s,$$

$$\overline{Y}_t = \{\xi^{(3)}(X_T) + E[\nabla_x \Phi(x, X_T)]_{x = X_T} \overline{X}_T + E[\nabla_{x'} \Phi(x, X_T) \overline{X}_T]_{x = X_T}\}$$

$$+ \int_t^T (\xi_s^{(4)}(\Lambda_s)$$

$$\text{(4.11)} \qquad\qquad + E[\nabla_\lambda f(\lambda, \Lambda_s)]_{\lambda = \Lambda_s} \overline{\Lambda_s} + E[\nabla_{\lambda'} f(\lambda, \Lambda_s) \overline{\Lambda_s}]_{\lambda = \Lambda_s}) \, ds$$

$$- \int_t^T \overline{Z}_s \, dW_s, \qquad t \in [0, T].$$

*Recall that $\overline{\mathbf{F}}$ is the filtration $\mathbf{F}$ augmented by the $\sigma$-filed generated by the random field $\xi$. Then the sequence $(\sqrt{N}(X^N - X, Y^N - Y, Z^N - Z))_{N \ge 1}$ converges in law over $C([0, T]; \mathbb{R}^d) \times C([0, T]; \mathbb{R}) \times L^2([0, T]; \mathbb{R}^d)$ to $\overline{\Lambda} = (\overline{X}, \overline{Y}, \overline{Z})$.*

REMARK 4.2.  As in the case of the only forward equation (cf. Remark 4.1) we get the pathwise continuity of the above introduced zero-mean Gaussian field with the help of Kolmogorov's Continuity Criterion for multi-parameter processes. From Remark 4.1 we know already that for all $m \ge 1$ there is some constant $C_m$ such that

$$E[|\xi_t^{(i)}(x) - \xi_{t'}^{(i)}(x')|^{2m}] \le C_m(|t - t'|^m + |x - x'|^{2m})$$

for all $(t, x), (t', x') \in [0, T] \times \mathbb{R}^d$, $i = 1, 2$. The same argument also shows that, again for some constant $C_m$,

$$E[|\xi^{(3)}(x) - \xi^{(3)}(x')|^{2m}] \le C_m |x - x'|^{2m} \qquad \text{for all } x, x' \in \mathbb{R}^d.$$

Moreover, for some generic constant $C_m$ we have

$$E[|\xi_t^{(4)}(x, y, z) - \xi_{t'}^{(4)}(x', y', z')|^{2m}]$$

$$\le C_m(E[|\xi_t^{(4)}(x, y, z) - \xi_{t'}^{(4)}(x', y', z')|^2])^m$$

$$= C_m(E[\text{Var}(f((x, y, z), \Lambda_t) - f((x', y', z'), \Lambda_{t'}))])^m$$

$$\leq C_m(|x - x'|^{2m} + |y - y'|^{2m} + |z - z'|^{2m}$$
$$+ E[|X_t - X_{t'}|^{2m} + |Y_t - Y_{t'}|^{2m}])$$

for all $t, t' \in [0, T], (x, y, z), (x', y', z') \in \mathbb{R}^d \times \mathbb{R} \times \mathbb{R}^d$ [recall the assumption (E)]. From standard SDE estimates we have

$$E[|X_t - X_{t'}|^{2m}] \leq C_m |t - t'|^m, \qquad t, t' \in [0, T].$$

Since, on the other hand, due to Corollary 4.1, we also have

$$E[|Y_t - Y_{t'}|^{2m}] \leq C_m |t - t'|^m \qquad \text{for all } t, t' \in [0, T],$$

we can conclude that

$$E[|\xi_t^{(4)}(x, y, z) - \xi_{t'}^{(4)}(x', y', z')|^{2m}]$$
$$\leq C_m(|t - t'|^m + |x - x'|^{2m} + |y - y'|^{2m} + |z - z'|^{2m})$$

for all $t, t' \in [0, T], (x, y, z), (x', y', z') \in \mathbb{R}^d \times \mathbb{R} \times \mathbb{R}^d$.

The proof of Theorem 4.2 is split in a sequel of statements which translate the proof of Theorem 4.1 into the context of a system composed of a forward and a backward SDE. For this we note that we can characterize $\sqrt{N}(Y^N - Y, Z^N - Z)$ as unique solution of the following BSDE:

$$\sqrt{N}(Y_t^N - Y_t)$$
$$= \xi^{3,N}(X_T^N) + \sqrt{N}(E[\Phi(x, X_T^N)]_{x=X_T^N} - E[\Phi(x, X_T)]_{x=X_T})$$
$$+ \int_t^T \{\xi_s^{4,N}(\Lambda_s^N) + \sqrt{N}(E[f(\lambda, \Lambda_s^N)]_{\lambda=\Lambda_s^N} - E[f(\lambda, \Lambda_s)]_{\lambda=\Lambda_s})\} \, ds$$
$$- \int_t^T \sqrt{N}(Z_s^N - Z_s) \, dW_s, \qquad t \in [0, T],$$

where

$$\xi^{3,N}(x) = \frac{1}{\sqrt{N}} \sum_{k=1}^N (\Phi(x, \Theta^k(X_T^N)) - E[\Phi(x, X_T^N)]),$$
$$(t, x) \in [0, T] \times \mathbb{R}^d,$$

$$\xi_t^{4,N}(\lambda) = \frac{1}{\sqrt{N}} \sum_{k=1}^N (f(\lambda, \Theta^k(\Lambda_t^N)) - E[f(\lambda, \Lambda_t^N)]),$$
$$(t, \lambda) \in [0, T] \times \mathbb{R}^d \times \mathbb{R} \times \mathbb{R}^d.$$

Recall also the definition of the random fields $\xi^{i,N} = \{\xi_t^{i,N}(x), (t, x) \in [0, T] \times \mathbb{R}^d\}$, $i = 1, 2$, which have been used for rewriting the (forward) equation for $X^N$. We remark that the random field $\xi^{i,N} = \{\xi_t^{i,N}(x), (t,x) \in [0,T] \times \mathbb{R}^d\}$, $i = 1, 2$, $\xi^{3,N} = \{\xi^{3,N}(x), x \in \mathbb{R}^d\}$ and $\xi^{4,N} = \{\xi_t^{4,N}(\lambda), (t,\lambda) \in [0,T] \times \mathbb{R}^d \times \mathbb{R} \times \mathbb{R}^d\}$ are independent of the Brownian motion $W$ and their paths are jointly continuous in all parameters. Moreover, we can characterize their limit behavior as follows:

PROPOSITION 4.3. *The sequence of the laws of the stochastic fields $\xi^N = (\xi^{1,N}, \xi^{2,N}, \xi^{3,N}, \xi^{4,N})$, $N \geq 1$, converges weakly on*

$$\mathcal{C} := C([0,T] \times \mathbb{R}^d; \mathbb{R}^d) \times C([0,T] \times \mathbb{R}^d; \mathbb{R}^{d \times d}) \times C(\mathbb{R}^d; \mathbb{R})$$
$$\times C([0,T] \times \mathbb{R}^d \times \mathbb{R} \times \mathbb{R}^d; \mathbb{R})$$

*to the law of the continuous zero-mean Gaussian field $\xi = (\xi^{(1)}, \xi^{(2)}, \xi^{(3)}, \xi^{(4)})$ introduced in Theorem 4.2.*

The proof of the proposition uses the same approach as that developed for the proof of Proposition 4.1. In a first lemma we establish the tightness of the sequence of laws $P_{\xi^N}$, $N \geq 1$, on $\mathcal{C}$, and a second lemma is devoted to the convergence of the finite-dimensional distributions of the sequence $\{\xi^N, N \geq 1\}$. Both statements together yield the weak convergence of the laws $P_{\xi^N}$, $N \geq 1$, on $\mathcal{C}$.

LEMMA 4.4. *The sequence of the laws of the stochastic processes $\xi^N = (\xi^{1,N}, \xi^{2,N}, \xi^{3,N}, \xi^{4,N})$, $N \geq 1$, is tight on $\mathcal{C}$. In fact, in addition to the estimates established in Lemma 4.1 for $\xi^{i,N}$, $N \geq 1, i = 1, 2$, we have for all $m \geq 1$ the existence of some constant $C_m$ such that:*

(i)          $E[|\xi^{3,N}(x)|^{2m} + |\xi_t^{4,N}(x,y,z)|^{2m}] \leq C_m$

*for all $(t,x,y,z) \in [0,T] \times \mathbb{R}^d \times \mathbb{R} \times \mathbb{R}^d, N \geq 1$;*

(ii)   $E[|\xi^{3,N}(x) - \xi^{3,N}(x')|^{2m} + |\xi_t^{4,N}(x,y,z) - \xi_{t'}^{4,N}(x',y',z')|^{2m}]$
$$\leq C_m(|t - t'|^m + |x - x'|^{2m} + |y - y'|^{2m} + |z - z'|^{2m}),$$

*for all $(t,x,y,z), (t',x',y',z') \in [0,T] \times \mathbb{R}^d \times \mathbb{R} \times \mathbb{R}^d, N \geq 1$.*

*Moreover, for all $m \geq 1$ there exists of some constant $C_m$ such that*

(iii)   $E\left[ \displaystyle\sup_{(x,y,z) \in \mathbb{R}^d \times \mathbb{R} \times \mathbb{R}^d} \left( \frac{|\xi^{3,N}(x)| + |\xi_t^{4,N}(x,y,z)|}{1 + |x| + |y| + |z|} \right)^{2m} \right] \leq C_m$

*for all $t \in [0,T], N \geq 1$.*

PROOF.    The proof of the estimates (i) and (iii) uses the same arguments as those developed in the proof of (i) and (iii) of Lemma 4.1, with the only difference that now Morrey's inequality has to be applied for $m \geq d + 1$ and $\gamma = 1 - \frac{2d+1}{2m}$. As concerns the proof of estimate (ii) we follow the argument developed for the proof of (ii) of Lemma 4.1. So we get, for a generic constant $C_m$ and all $(t,x,y,z), (t',x',y',z') \in [0,T] \times \mathbb{R}^d \times \mathbb{R} \times \mathbb{R}^d$ and $N \geq 1$,

$$E[|\xi^{3,N}(x) - \xi^{3,N}(x')|^{2m} + |\xi_t^{4,N}(x,y,z) - \xi_{t'}^{4,N}(x',y',z')|^{2m}]$$
$$\leq C_m(E[|\Phi(x,X_T^N) - \Phi(x',X_T^N)|^{2m}]$$
$$+ E[|f((x,y,z),\Lambda_t^N) - f((x',y',z'),\Lambda_{t'}^N)|^{2m}])$$
$$\leq C_m(|x-x'|^{2m} + |y-y'|^{2m} + |z-z'|^{2m}$$
$$+ E[|X_t^N - X_{t'}^N|^{2m} + |Y_t^N - Y_{t'}^N|^{2m}]).$$

Recall that the coefficients $\Phi$ and $f$ are Lipschitz continuous and the assumption (E). From Corollary 4.1 and a standard estimate for forward SDEs we then get (ii).

Then, due to Kolmogorov's weak compactness criterion for multi-parameter processes the estimates (i) and (ii) of the present lemma and those of Lemma 4.1 imply the tightness of the sequence of laws of $\xi^N$, $N \geq 1$, on $\mathcal{C}$. $\quad\square$

For proving Proposition 4.3 we have still to show the following.

LEMMA 4.5. *The finite-dimensional laws of the stochastic fields $\xi^N = (\xi^{1,N}, \xi^{2,N}, \xi^{3,N}, \xi^{4,N})$, $N \geq 1$, converge weakly to the corresponding finite-dimensional laws of the continuous zero-mean Gaussian multi-parameter process $\xi = (\xi^{(1)}, \xi^{(2)}, \xi^{(3)}, \xi^{(4)})$ introduced in Theorem 4.2.*

PROOF. Following the idea of the proof of Lemma 4.2 we decompose $\xi^N = (\xi^{1,N}, \xi^{2,N}, \xi^{3,N}, \xi^{4,N})$ into the sum of the random fields $\xi^{\cdot,N,1} = (\xi^{1,N,1}, \xi^{2,N,1}, \xi^{3,N,1}, \xi^{4,N,1})$ and $\xi^{\cdot,N,2} = (\xi^{1,N,2}, \xi^{2,N,2}, \xi^{3,N,2}, \xi^{4,N,2})$, where the fields $\xi^{2,N,1}$ and $\xi^{2,N,2}$ have been introduced in the proof of Lemma 4.2, the fields $\xi^{1,N,1}$ and $\xi^{1,N,2}$ are defined in the same way, with $b$ instead of $\sigma$, and ditto for $\xi^{3,N,1}$ and $\xi^{3,N,2}$ where we have $\Phi$ instead of $\sigma$ and $t = T$. Finally, in the same spirit we define

$$\xi_t^{4,N,1}(\lambda) = \frac{1}{\sqrt{N}} \sum_{k=1}^{N} (\Theta^k(f(\lambda,\Lambda_t)) - E[f(\lambda,\Lambda_t)]),$$
$$(t,\lambda) \in [0,T] \times \mathbb{R}^d \times \mathbb{R} \times \mathbb{R}^d;$$

$$\xi_t^{4,N,2}(x) = \frac{1}{\sqrt{N}} \sum_{k=1}^{N} (\Theta^k(f(\lambda,\Lambda_t^N) - f(\lambda,\Lambda_t)) - E[f(\lambda,\Lambda_t^N) - f(\lambda,\Lambda_t)])$$
$$\text{for } (t,\lambda) \in [0,T] \times \mathbb{R}^d \times \mathbb{R} \times \mathbb{R}^d.$$

Since the bounded multi-parameter fields $\Theta^k(b(\cdot,X), \sigma(\cdot,X), \Phi(\cdot,X_T), f(\cdot, \Lambda))$, $k \geq 0$ are i.i.d., it follows directly from the central limit theorem that the finite-dimensional laws of the random field $\xi^{\cdot,N,1} = (\xi^{1,N,1}, \xi^{2,N,1}, \xi^{3,N,1}, \xi^{4,N,1})$ converge weakly to the corresponding finite-dimensional laws of the zero-mean Gaussian multi-parameter process $\xi = (\xi^{(1)}, \xi^{(2)}, \xi^{(3)}, \xi^{(4)})$ whose covariance functions coincides with those of the field $(b(\cdot, X), \sigma(\cdot, X), \Phi(\cdot, X_T), f(\cdot, \Lambda))$ [recall the definition of $\xi = (\xi^{(1)}, \xi^{(2)}, \xi^{(3)}, \xi^{(4)})$ given in Theorem 4.2]. On the other hand, by repeating the argument developed in the proof of Lemma 4.2 and recalling the assumption (E) we get from Proposition 3.2 and Theorem 3.2 that

$$E[|\xi_t^{1,N,2}(x)|^2] + E[|\xi_t^{2,N,2}(x)|^2] + E[|\xi^{3,N,2}(x)|^2] + E[|\xi_t^{4,N,2}(x,y,z)|^2]$$

$$= \sum_{i=1}^d \mathrm{Var}(b^i(x, X_t^N) - b^i(x, X_t)) + \sum_{i,j=1}^d \mathrm{Var}(\sigma^{i,j}(x, X_t^N) - \sigma^{i,j}(x, X_t))$$

$$\quad + \mathrm{Var}(\Phi(x, X_T^N) - \Phi(x, X_T))$$

$$\quad + \mathrm{Var}(f((x,y,z), \Lambda_t^N) - f((x,y,z), \Lambda_t))$$

$$\leq C(E[|X_t^N - X_t|^2] + E[|X_T^N - X_T|^2] + E[|Y_t^N - Y_t|^2])$$

$$\leq \frac{C}{N} \qquad \text{for all } N \geq 1, (t, x, y, z) \in [0, T] \times \mathbb{R}^d \times \mathbb{R} \times \mathbb{R}^d.$$

Consequently, it follows that also the finite-dimensional laws of the random fields $\xi^N = (\xi^{1,N}, \xi^{2,N}, \xi^{3,N}, \xi^{4,N})$, $N \geq 1$, converge weakly to the corresponding finite-dimensional laws of the zero-mean Gaussian field $\xi = (\xi^{(1)}, \xi^{(2)}, \xi^{(3)}, \xi^{(4)})$. The proof is complete. $\square$

As we have already pointed out above, Proposition 4.3 follows directly from the Lemmas 4.4 and 4.5. On the other hand, Proposition 4.3 allows to adapt the approach for the study of the limit behavior of $\sqrt{N}(X^N - X)$ to that of the triplet $\sqrt{N}(X^N - X, Y^N - Y, Z^N - Z)$: using the same notation as in the preceding section we remark that, since the random fields $\xi^N$, $N \geq 1$, are independent of the driving Brownian motion $W$, Skorohod's representation theorem yields that on an appropriate complete probability space $(\Omega', \mathcal{F}', P')$ there exist copies $\xi'^N$, $N \geq 1$, and $\xi'$ of the stochastic fields $\xi^N$, $N \geq 1$, and $\xi$, as well as a $d$-dimensional Brownian motion $W'$, such that:

 (i) $P'_{\xi'} = P_\xi, P'_{\xi'^N} = P_{\xi^N}, N \geq 1$;

 (ii) $\xi'^N = (\xi'^{1,N}, \xi'^{2,N}, \xi'^{3,N}, \xi'^{4,N}) \longrightarrow \xi' = (\xi'^1, \xi'^2, \xi'^3, \xi'^4)$, uniformly on the compacts of $[0, T] \times \mathbb{R}^d \times \mathbb{R} \times \mathbb{R}^d$, $P'$-a.s.;

 (iii) $W'$ is independent of $\xi'$ and $\xi'^N$, for all $N \geq 1$.

The probability space $(\Omega', \mathcal{F}', P')$ is endowed with the filtration $\mathbf{F}' = (\mathcal{F}'_t = \mathcal{F}_t^{W'} \vee \mathcal{F}'_0)_{t \in [0,T]}$, where $\mathcal{F}'_0 = \sigma\{\xi'_s(x), \zeta'^N_s(x), (s, x) \in [0, T] \times \mathbb{R}, N \geq 1\} \vee \mathcal{N}$.

Remark that, with respect to this filtration, the process $W'$ is a Brownian motion and has the martingale representation property.

Given $\xi'^N, \xi'$ and the Brownian motion $W'$ we redefine the triplets $(X^N, Y^N, Z^N)$ and $(X, Y, Z)$ on the new probability space: the processes $X^N$ and $X$ are redefined as in the preceding section (their redefinitions are denoted again by $X'^N$ and $X'$, resp.), and in the same spirit we introduce the couples $(Y'^N, Z'^N)$, $(Y', Z') \in \mathcal{S}_{\mathbf{F}'}^2([0,T];\mathbb{R}) \times L_{\mathbf{F}'}^2([0,T];\mathbb{R}^d)$ as unique solution of the backward equations

$$Y_t'^N = \left( \frac{1}{\sqrt{N}} \xi'^{3,N}(X_T'^N) + E'[\Phi(x, X_T'^N)]_{|x=X_T'^N} \right)$$
$$+ \int_t^T \left( \frac{1}{\sqrt{N}} \xi_s'^{4,N}(\Lambda_s'^N) + E'[f(\lambda, \Lambda_s'^N)]_{|\lambda=\Lambda_s'^N} \right) ds - \int_t^T Z_s'^N dW_s',$$
$$t \in [0,T]$$

and

$$Y_t' = E'[\Phi(x, X_T')]_{|x=X_T'} + \int_t^T E'[f(\lambda, \Lambda_s')]_{|\lambda=\Lambda_s'} ds - \int_t^T Z_s' dW_s',$$
$$t \in [0,T],$$

respectively, where $\Lambda'^N = (X'^N, Y'^N, Z'^N)$ and $\Lambda' = (X', Y', Z')$. Recall that the coefficients $\xi^{3,N}(\cdot)$ and $\xi_s^{4,N}(\cdot)$ are Lipschitz, uniformly with respect to $(\omega, s) \in \Omega \times [0,T]$, and so are $\xi'^{3,N}(\cdot)$ and $\xi_s'^{4,N}(\cdot)$ with respect to $(\omega', s) \in \Omega' \times [0,T]$. Thus, since $P' \circ (\xi'^N, W')^{-1} = P \circ (\xi^N, W)^{-1}$, $N \geq 1$, and $P' \circ (\xi', W')^{-1} = P \circ (\xi, W)^{-1}$, we can conclude that also $P' \circ (\xi'^N, W', \Lambda'^N)^{-1} = P \circ (\xi^N, W, \Lambda^N)^{-1}$, for all $N \geq 1$, and $P' \circ (\xi', W', \Lambda')^{-1} = P \circ (\xi, W, \Lambda)^{-1}$. This allows to reduce the study of the limit behavior of $\sqrt{N}(\Lambda^N - \Lambda), N \geq 1$, to that of the sequence $\sqrt{N}(\Lambda'^N - \Lambda'), N \geq 1$.

PROPOSITION 4.4. *Let* $\widehat{X} \in \mathcal{S}_{\mathbf{F}'}^2([0,T];\mathbb{R}^d)$ *be the unique solution of SDE* (4.5) *and* $(\widehat{Y}, \widehat{Z}) \in \mathcal{S}_{\mathbf{F}'}^2([0,T];\mathbb{R}) \times L_{\mathbf{F}'}^2([0,T];\mathbb{R}^d)$ *that of the backward equation*

$$\widehat{Y}_t = \{\xi'^3(X_T') + E'[\nabla_x \Phi(x, X_T')]_{x=X_T'} \widehat{X}_T$$
$$+ E'[\nabla_{x'} \Phi(x, X_T') \widehat{X}_T]_{x=X_T'}\}$$
(4.12)
$$+ \int_t^T (\xi_s'^4(\Lambda_s') + E'[\nabla_\lambda f(\lambda, \Lambda_s')]_{\lambda=\Lambda_s'} \widehat{\Lambda}_s$$
$$+ E[\nabla_{\lambda'} f(\lambda, \Lambda_s') \widehat{\Lambda}_s]_{\lambda=\Lambda_s'}) ds$$
$$- \int_t^T \widehat{Z}_s dW_s, \qquad t \in [0,T], \text{ where } \widehat{\Lambda} = (\widehat{X}, \widehat{Y}, \widehat{Z}).$$

*Then*

$$E'\Big[\sup_{s\in[0,T]}(|\sqrt{N}(X_s'^N - X_s') - \widehat{X}_s|^2 + |\sqrt{N}(Y_s'^N - Y_s') - \widehat{Y}_s|^2)$$

$$+ \int_0^T |\sqrt{N}(Z_s'^N - Z_s') - \widehat{Z}_s|^2\, ds\Big] \longrightarrow 0$$

$$as\ N \to +\infty.$$

PROOF. We begin by remarking that an argument similar to that of Step 1 of the proof of Proposition 4.2, combined with estimate (i) of Lemma 4.4, allows to show that $\xi'^3(X_T') \in L^p(\Omega', \mathcal{F}_T', P')$ and $\xi'^4(\Lambda') \in L^p_{\mathbf{F}'}([0,T]; \mathbb{R})$, for all $p \geq 1$. Consequently, since the coefficients $\nabla_x\Phi$, $\nabla_{x'}\Phi$, $\nabla_\lambda f$ and $\nabla_{\lambda'} f$ are bounded, the above backward SDE admits a unique solution $(\widehat{Y}, \widehat{Z})$. Moreover, this solution belongs even to $\mathcal{S}^p_{\mathbf{F}'}([0,T]; \mathbb{R}) \times L^p_{\mathbf{F}'}(\Omega'; L^2([0,T]; \mathbb{R}^d))$, for all $p \geq 1$.

From Step 2 of the proof of Proposition 4.2 we know already that $(\xi_t^{i,N}(X_t'^N))_{t\in[0,T]} \to (\xi_t'^{(i)}(X_t'))_{t\in[0,T]}$ in $L^2([0,T] \times \Omega', dt\, dP')$, as $N \to +\infty$. The same argument yields that also $(\xi'^{3,N}(X_T'^N)) \to (\xi'^{(3)}(X_T'))$ in $L^2(\Omega', \mathcal{F}', P')$, as $N \to +\infty$. Moreover, $s \in [0,T]$,

$$E'[|\xi_s'^{4,N}(\Lambda_s'^N)|^m]$$
$$= E[|\xi_s^{4,N}(\Lambda_s^N)|^m]$$
$$\leq \left(E\Big[\sup_{\lambda\in\mathbb{R}^d\times\mathbb{R}\times\mathbb{R}^d}\Big(\frac{|\xi_s^{4,N}(\lambda)|}{1+|\lambda|}\Big)^{2m}\Big]\right)^{1/2} \left(E[(1+|\Lambda_s^N|)^{2m}]\right)^{1/2}.$$

Thus, using Lemma 4.4(iii) and the fact that, for all $N, m \geq 1$,

$$E[|\Lambda_s^N|^{2m}] \leq C_m E[|X_s^N|^{2m} + |Y_s^N|^{2m} + |Z_s^N|^{2m}] \leq C_m, \qquad ds\text{-a.e.}$$

Recall that the coefficients of (4.1) and of BSDE (4.7) are bounded, uniformly with respect to $N \geq 1$; the same holds true for $Z^N$ (cf. Lemma 4.3), we can conclude that, for some constant $C_m$,

$$E'[|\xi_s'^{4,N}(\Lambda_s'^N)|^m] \leq C_m, \qquad ds\text{-a.e. on } [0,T], N \geq 1.$$

On the other hand, recalling that $\xi'^N = (\xi'^{1,N}, \xi'^{2,N}, \xi'^{3,N}, \xi'^{4,N}) \to \xi' = (\xi'^1, \xi'^2, \xi'^3, \xi'^4)$ uniformly on the compacts of $[0,T] \times \mathbb{R}^d \times \mathbb{R} \times \mathbb{R}^d$, $P'$-a.s., and $N \geq 1$,

$$E'\Big[\sup_{s\in[0,T]}(|X_s'^N - X_s'|^2 + |Y_s'^N - Y_s'|^2) + \int_0^T |Z_t'^N - Z_t'|^2\, dt\Big]$$

$$= E\Big[\sup_{s\in[0,T]}(|X_s^N - X_s|^2 + |Y_s^N - Y_s|^2) + \int_0^T |Z_t^N - Z_t|^2\, dt\Big] \leq \frac{C}{N}$$

(cf. Proposition 3.2 and Theorem 3.2), we see that

$$\sum_{i=1}^{2} \sup_{s \in [0,T]} |\xi_s'^{i,N}(X_s'^N) - \xi_s'^i(X_s')| + |\xi'^{3,N}(X_T'^N) - \xi'^3(X_T')| \longrightarrow 0$$

in probability $P'$, and $\xi_s'^{4,N}(\Lambda_s'^N) \longrightarrow \xi_s'^4(\Lambda_s')$, in measure $ds\,dP'$.

Combining this with the uniform square integrability of $\{\xi_s'^N(X_s'^N), s \in [0,T], N \geq 1\}$, proved above, yields the convergence of $\xi'^{i,N}(X'^N)$ to $\xi'^i(X')$ in $L^2([0,T] \times \Omega', dt\,dP')$, for $i = 1, 2$, that of $\xi'^{3,N}(X'^N)$ to $\xi'^3(X')$ in $L^2(\Omega', \mathcal{F}', P')$ and that of $\xi'^{4,N}(\Lambda'^N)$ to $\xi'^4(\Lambda')$ in $L^2([0,T] \times \Omega', dt\,dP')$.

Let us now prove the convergence of $\sqrt{N}(X'^N - X', Y'^N - Y', Z'^N - Z')$ to $\widehat{\Lambda} = (\widehat{X}, \widehat{Y}, \widehat{Z})$. From Proposition 4.2 we know already that

$$E'\left[ \sup_{s \in [0,t]} |\sqrt{N}(X_s'^N - X_s') - \widehat{X}_s|^2 \right] \longrightarrow 0 \qquad \text{as } N \to +\infty.$$

For verifying the convergence of $\sqrt{N}(Y'^N - Y', Z'^N - Z')$ we observe that

$$\sqrt{N}(Y_t'^N - Y_t') - \widehat{Y}_t = I_t^{1,N} + I^{2,N} + I_t^{3,N} - \int_t^T \{\sqrt{N}(Z_s'^N - Z_s') - \widehat{Z}_s\}\,dW_s',$$
$$t \in [0,T],$$

where

$$I_t^{1,N} = \{\xi'^{3,N}(X_T'^N) - \xi'^3(X_T')\} + \int_t^T (\xi_s'^{4,N}(\Lambda_s'^N) - \xi'^4(\Lambda_s'))\,ds,$$

$$I^{2,N} = \sqrt{N}(E'[\Phi(x, X_T'^N)]_{x=X_T'^N} - E'[\Phi(x, X_T')]_{x=X_T'})$$
$$\qquad - (E'[\nabla_x \Phi(x, X_T')]_{x=X_T'} \widehat{X}_T + E'[\nabla_{x'} \Phi(x, X_T') \widehat{X}_T]_{x=X_T'}),$$

$$I_t^{3,N} = \int_t^T \{\sqrt{N}(E'[f(\lambda, \Lambda_s'^N)]_{\lambda=\Lambda_s'^N} - E'[f(\lambda, \Lambda_s')]_{\lambda=\Lambda_s'})$$
$$\qquad - (E'[\nabla_\lambda f(\lambda, \Lambda_s')]_{\lambda=\Lambda_s'} \widehat{\Lambda}_s + E'[\nabla_{\lambda'} f(\lambda, \Lambda_s') \widehat{\Lambda}_s]_{\lambda=\Lambda_s'})\}\,ds,$$
$$t \in [0,T].$$

From the convergence of $\xi'^{3,N}(X'^N)$ to $\xi'^3(X')$ in $L^2(\Omega', \mathcal{F}', P')$ and that of $\xi'^{4,N}(\Lambda'^N)$ to $\xi'^4(\Lambda')$ in $L^2([0,T] \times \Omega', dt\,dP')$ we obtain that

$$E'\left[ \sup_{t \in [0,T]} |I_t^{1,N}|^2 \right] \longrightarrow 0 \qquad \text{as } N \to +\infty.$$

In analogy to Step 3 of the proof of Proposition 4.2 we get that, for some real sequence $0 < \rho_N \searrow 0$ $(N \to +\infty)$,

$$E'[|I^{2,N}|^2] \leq \rho_N + CE'\left[ \left| \sqrt{N}(X_T'^N - X_T') - \widehat{X}_T \right|^2 \right]$$

$$\leq \rho_N + CE'\Big[\sup_{t\in[0,T]}|\sqrt{N}(X_t'^N - X_t') - \widehat{X}_t|^2\Big] \longrightarrow 0$$

$$\text{as } N \to +\infty.$$

It remains to estimate $I^{3,N}$. For this, thanks to the continuous differentiability of $f$, with the notation $\Lambda_s'^{N,\gamma} := \Lambda_s' + \gamma(\Lambda_s'^N - \Lambda_s')$ we get

$$E'\Big[\sup_{s\in[t,T]}|I_s^{3,N}|^2\Big]$$

$$\leq CE'\Big[\int_t^T \Big|\int_0^1 E'[\nabla_\lambda f(\lambda_1, \Lambda_s'^{N,\gamma})$$

$$- \nabla_\lambda f(\lambda_2, \Lambda_s')]_{|\lambda_1 = \Lambda_s'^{N,\gamma}, \lambda_2 = \Lambda_s'}\, d\gamma \times \widehat{\Lambda}_s$$

$$+ \int_0^1 E'[\{\nabla_{\lambda'} f(\lambda_1, \Lambda_s'^{N,\gamma})$$

$$- \nabla_{\lambda'} f(\lambda_2, \Lambda_s')\}\widehat{\Lambda}_s]_{|\lambda_1 = \Lambda_s'^{N,\gamma}, \lambda_2 = \Lambda_s'}\, d\gamma\Big|^2\, ds\Big]$$

$$+ C(T-t)E'$$

$$\times \Big[\int_t^T |E'[\nabla_\lambda f(\lambda, \Lambda_s'^{N,\gamma})]_{\lambda = \Lambda_s'}(\sqrt{N}(\Lambda_s'^N - \Lambda_s') - \widehat{\Lambda}_s)$$

$$+ E'[\nabla_{\lambda'} f(\lambda, \Lambda_s'^{N,\gamma})(\sqrt{N}(\Lambda_s'^N - \Lambda_s') - \widehat{\Lambda}_s)]_{\lambda = \Lambda_s'}|^2\, ds\Big],$$

$$t\in[0,T],$$

and taking into account that the first-order derivatives of $f$ are bounded and continuous, it follows from the convergence of $\Lambda'^{N,\gamma} \to \Lambda'$ $(N \to +\infty)$ in $\mathcal{B}^2_{\mathbf{F}'}$, for all $\gamma \in [0,1]$, that for some sequence $0 < \rho_N \searrow 0$,

$$E'\Big[\sup_{s\in[t,T]}|I_s^{3,N}|^2\Big] \leq \rho_N + C(T-t)\int_t^T E'[|\sqrt{N}(\Lambda_s'^N - \Lambda_s') - \widehat{\Lambda}_s|^2]\, ds,$$

$$t\in[0,T], N\geq 1.$$

Thus, for some sufficiently small $\delta > 0$ which depends only on the Lipschitz constant of $f$, and for some sequence $0 < \rho_N \searrow 0$, we have for all $t \in [T - \delta, T]$,

$$E'[|\sqrt{N}(Y_t'^N - Y_t') - \widehat{Y}_t|^2] + \frac{1}{2}E'\Big[\int_t^T |\sqrt{N}(Z_s'^N - Z_s') - \widehat{Z}_s|^2\, ds\Big]$$

$$\leq \rho_N + CE'\Big[\int_t^T |\sqrt{N}(X_s'^N - X_s') - \widehat{X}_s|^2\, ds\Big]$$

$$+ CE'\left[\int_t^T |\sqrt{N}(Y_s'^N - Y_s') - \widehat{Y}_s|^2\, ds\right]$$

and from Proposition 4.2 and Gronwall's inequality we get

$$E'[|\sqrt{N}(Y_t'^N - Y_t') - \widehat{Y}_t|^2] \to 0 \qquad \text{for all } t \in [T - \delta, T]$$

and

$$E'\left[\int_{T-\delta}^T |\sqrt{N}(Z_s'^N - Z_s') - \widehat{Z}_s|^2\, ds\right] \to 0 \qquad \text{as } N \to +\infty.$$

Having got the convergence of $\sqrt{N}(Y_{T-\delta}'^N - Y_{T-\delta}') - \widehat{Y}_{T-\delta}$ in $L^2(\Omega', \mathcal{F}', P')$ we can repeat the above convergence argument on $[T - 2\delta, T - \delta]$ by replacing $I^{2,N}$ by $\sqrt{N}(Y_{T-\delta}'^N - Y_{T-\delta}') - \widehat{Y}_{T-\delta}$. This yields $E'[|\sqrt{N}(Y_t'^N - Y_t') - \widehat{Y}_t|^2] \to 0$, for all $t \in [T - 2\delta, T - \delta]$, and $E'[\int_{T-2\delta}^{T-\delta} |\sqrt{N}(Z_s'^N - Z_s') - \widehat{Z}_s|^2\, ds] \to 0$. Iterating this argument we get the convergence over the whole interval $[0, T]$. Finally, from the equation

$$\sqrt{N}(Y_t'^N - Y_t') - \widehat{Y}_t = I_t^{1,N} + I^{2,N} + I_t^{3,N} - \int_t^T \{\sqrt{N}(Z_s'^N - Z_s') - \widehat{Z}_s\}\, dW_s',$$

$t \in [0, T]$, we get

$$E'\left[\sup_{t \in [0,T]} |\sqrt{N}(Y_t'^N - Y_t') - \widehat{Y}_t|^2\right]$$

$$\leq 4\left(E'\left[\sup_{t \in [0,T]} |I_t^{1,N}|^2\right] + E'[|I^{2,N}|^2]\right.$$

$$\left. + E'\left[\sup_{t \in [0,T]} |I_t^{3,N}|^2\right] + E'\left[\int_0^T |\sqrt{N}(Z_t'^N - Z_t') - \widehat{Z}_t|^2\, dt\right]\right)$$

$$\leq \rho_N + CE'\left[\int_0^T |\sqrt{N}(\Lambda_s'^N - \Lambda_s') - \widehat{\Lambda}_s|^2\, ds\right] \longrightarrow 0 \qquad \text{as } N \to +\infty.$$

The proof is complete. $\square$

Finally, it remains to give the proof of Theorem 4.2.

PROOF OF THEOREM 4.2. In analogy to the existence and uniqueness of the triplet $\widehat{\Lambda} = (\widehat{X}, \widehat{Y}, \widehat{Z}) \in \mathcal{B}_{\mathbb{F}'}^2$ we can prove that of $\overline{\Lambda} = (\overline{X}, \overline{Y}, \overline{Z}) \in \mathcal{B}_{\mathbb{F}}^2$. Since $W$ is independent of $\xi$ as well as $\xi^N, N \geq 1$, and, on the other hand, also $W'$ is independent of $\xi'$ and $\xi'^N, N \geq 1$, it follows that:

(i) $P \circ (W, \xi^N)^{-1} = P' \circ (W', \xi'^N)^{-1}$, for all $N \geq 1$;

(ii) $P \circ (W, \xi)^{-1} = P' \circ (W', \xi')^{-1}$.

We notice that $\Lambda'^N$ is the unique solution of an forward–backward SDE governed by the $W'$ and $\xi'^N$, and $\Lambda^N$ is the unique solution of the same system of equations but driven by $W$ and $\xi^N$ instead of $W'$ and $\xi'^N$. On the other hand, $\Lambda'$ is the unique solution of the system (4.8)–(4.9) with $W'$ at the place of $W$. Thus,

(iii) $P \circ (W, \xi^N, \Lambda, \Lambda^N)^{-1} = P' \circ (W', \xi'^N, \Lambda', \Lambda'^N)^{-1}$, for all $N \geq 1$.

Concerning the process $\widetilde{\Lambda}$, it is the unique solution of the same equation as $\overline{\Lambda}$, but with $\xi$ and $W$ replaced by $\xi'$ and $W'$. This, together with (ii), yields:

(iv) $P \circ (W, \xi, \Lambda, \overline{\Lambda})^{-1} = P' \circ (W', \xi', \Lambda', \widehat{\Lambda})^{-1}$.

From (iii) and (iv) we see, in particular, that $P_{\sqrt{N}(\Lambda^N - \Lambda)} = P'_{\sqrt{N}(\Lambda'^N - \Lambda')}$, $N \geq 1$, and $P_{\overline{\Lambda}} = P'_{\widehat{\Lambda}}$. For completing the proof it suffices now to remark that the convergence of $\sqrt{N}(\Lambda'^N - \Lambda')$ to $\widehat{\Lambda}$ in $\mathcal{S}^2_{\mathbf{F}'}([0, T]; \mathbb{R}^d)$ implies that of their laws on $C([0, T]; \mathbb{R}^d) \times C([0, T]; \mathbb{R}) \times L^2([0, T]; \mathbb{R}^d)$. The proof is complete. $\square$

## REFERENCES

[1] BOSSY, M. (2005). Some stochastic particle methods for nonlinear parabolic PDEs. In *GRIP—Research Group on Particle Interactions. ESAIM Proceedings* **15** 18–57. MR2441320

[2] BOSSY, M. and TALAY, D. (1997). A stochastic particle method for the McKean–Vlasov and the Burgers equation. *Math. Comp.* **66** 157–192. MR1370849

[3] BUCKDAHN, R., LI, J. and PENG, S. (2007). Mean-field backward stochastic differential equations and related patial differential equations. Submitted. Available at http://arxiv.org/abs/0711.2167.

[4] CHAN, T. (1994). Dynamics of the McKean–Vlasov equation. *Ann. Probab.* **22** 431–441. MR1258884

[5] EL KAROUI, N., PENG, S. and QUENEZ, M. C. (1997). Backward stochastic differential equations in finance. *Math. Finance* **7** 1–71. MR1434407

[6] KOTELENEZ, P. (1995). A class of quasilinear stochastic partial differential equations of McKean–Vlasov type with mass conservation. *Probab. Theory Related Fields* **102** 159–188. MR1337250

[7] LASRY, J.-M. and LIONS, P.-L. (2007). Mean field games. *Japan J. Math.* **2** 229–260. MR2295621

[8] MÉLÉARD, S. (1996). Asymptotic behaviour of some interacting particle systems; McKean–Vlasov and Boltzmann models. In *Probabilistic Models for Nonlinear Partial Differential Equations (Montecatini Terme, 1995). Lecture Notes in Math.* **1627** 42–95. Springer, Berlin. MR1431299

[9] NUALART, D. (1995). *The Malliavin Calculus and Related Topics.* Springer, New York. MR1344217

[10] OVERBECK, L. (1995). Superprocesses and McKean–Vlasov equations with creation of mass. Preprint.

[11] PARDOUX, É. and PENG, S. G. (1990). Adapted solution of a backward stochastic differential equation. *Systems Control Lett.* **14** 61–74. MR1037747

[12] PARDOUX, É. and PENG, S. (1992). Backward stochastic differential equations and quasilinear parabolic partial differential equations. In *Stochastic Differential Equations and Their Applications (Charlotte, NC, 1991). Lecture Notes in Control and Information Sciences* **176** 200–217. Springer, Berlin. MR1176785

[13] PENG, S. G. (1992). A generalized dynamic programming principle and Hamilton–Jacobi–Bellman equation. *Stochastics Rep.* **38** 119–134. MR1274898

[14] PRA, P. D. and HOLLANDER, F. D. (1995). McKean–Vlasov limit for interacting random processes in random media. *J. Statist. Phys.* **84** 735–772. MR1400186

[15] SZNITMAN, A.-S. (1984). Nonlinear reflecting diffusion process, and the propagation of chaos and fluctuations associated. *J. Funct. Anal.* **56** 311–336. MR743844

[16] SZNITMAN, A.-S. (1991). Topics in propagation of chaos. In *École D'Été de Probabilités de Saint-Flour XIX—1989. Lecture Notes in Math.* **1464** 165–252. Springer, Berlin. MR1108185

[17] TALAY, D. and VAILLANT, O. (2003). A stochastic particle method with random weights for the computation of statistical solutions of McKean–Vlasov equations. *Ann. Appl. Probab.* **13** 140–180. MR1951996

R. BUCKDAHN
DÉPARTEMENT DE MATHÉMATIQUES
UNIVERSITÉ DE BRETAGNE OCCIDENTALE
6, AVENUE VICTOR-LE-GORGEU
B.P. 809, 29285 BREST CEDEX
FRANCE
E-MAIL: Rainer.Buckdahn@univ-brest.fr

B. DJEHICHE
DEPARTMENT OF MATHEMATICS
ROYAL INSTITUTE OF TECHNOLOGY
S-100 44 STOCKHOLM
SWEDEN
E-MAIL: boualem@math.kth.se

J. LI
SCHOOL OF MATHEMATICS AND STATISTICS
SHANDONG UNIVERSITY AT WEIHAI
180 WENHUA XILU, WEIHAI 264209
P. R. CHINA
E-MAIL: juanli@sdu.edu.cn

S. PENG
SCHOOL OF MATHEMATICS
AND SYSTEM SCIENCES
SHANDONG UNIVERSITY
JINAN 250100
P. R. CHINA
E-MAIL: peng@sdu.edu.cn